\documentclass[11pt]{article}

\usepackage{amsmath,amsxtra,amsthm,amssymb,makeidx,graphics}
\usepackage{latexsym,color}
\theoremstyle{plain}
\newtheorem{theorem}{Theorem}[section]
\newtheorem{lemma}[theorem]{Lemma}
\newtheorem{remark}[theorem]{Remark}
\newtheorem{corollary}[theorem]{Corollary}
\newtheorem{definition}[theorem]{Definition}
\newtheorem{example}[theorem]{Example}
\newtheorem{exercise}[theorem]{Exercise}

\theoremstyle{remark}

\newcommand{\ncom}{\newcommand}

\newcommand{\tab}{\mbox{tab}}

\newcommand{\cB}{{\cal B}}

\newcommand{\cD}{{\cal D}}

\newcommand{\cM}{{\cal M}}

\newcommand{\cO}{{\cal O}}
\newcommand{\cP}{{\cal P}}
\newcommand{\cbP}{{{\ol{\cal P}}}}
\newcommand{\cbY}{{{\ol{\cal Y}}}}

\newcommand{\cS}{{\cal S}}

\newcommand{\cY}{{\cal Y}}

\newcommand{\C}{{\mathbb C}}
\newcommand{\Q}{{\mathbb Q}}

\newcommand{\Z}{{\mathbb Z}}

\newcommand{\End}{{\mbox{End}}}

\newcommand{\ind}{{\mbox{ind}}}
\newcommand{\res}{{\mbox{res}}}
\newcommand{\supp}{{\mbox{supp}}}

\ncom{\n}[1]{{{ \|[#1]\|}}}
\ncom{\ns}{\normalsize}
\ncom{\la}{\lambda}
\ncom{\bm}{\boldmath}
\ncom{\noi}{\noindent}
\ncom{\bq}{\begin{equation}}
\ncom{\eq}{\end{equation}}  
\ncom{\beqn}{\begin{eqnarray*}}
\ncom{\eeqn}{\end{eqnarray*}}  
\ncom{\ba}{\begin{array}}
\ncom{\ul}{\underline}   

\ncom{\ea}{\end{array}}
\ncom{\beq}{\begin{eqnarray}}
\ncom{\eeq}{\end{eqnarray}}  
\ncom{\nno}{\nonumber}
\ncom{\hs}{\mbox{\hspace{.25cm}}}
\ncom{\rar}{\rightarrow}
\ncom{\lrar}{\longrightarrow}
\ncom{\Rar}{\Rightarrow}
\ncom{\noin}{\noindent} 
\ncom{\bc}{\begin{center}}
\ncom{\ec}{\end{center}}  
\ncom{\sz}{\scriptsize}   
\ncom{\fpd}{\Phi(\pi^{'})}
\ncom{\fp}{\Phi(\pi) }
\ncom{\nk}{\left< \begin{array}{c}
                       n\\k \end{array} \right>}
\ncom{\nd}{1^{'},2^{'},\cdots,n^{'}}
\ncom{\de}{\bigtriangleup (F_{2n},\leq)}
\ncom{\del}{\bigtriangleup}
\ncom{\cov}{<\!\!\!\!\cdot }
\ncom{\bt}{\begin{theorem}}
\ncom{\bcon}{\begin{con}}
\ncom{\et}{\end{theorem}}
\ncom{\econ}{\end{con}}
\ncom{\bl}{\begin{lemma}}
\ncom{\el}{\end{lemma}}  
\ncom{\bco}{\begin{corollary}} 
\ncom{\ds}{\displaystyle}
\ncom{\eco}{\end{corollary}}   
\ncom{\bp}{\begin{pro}}  
\ncom{\ep}{\end{pro}}    
\ncom{\bex}{\begin{example}}
\ncom{\eex}{\end{example}}  
\ncom{\bexr}{\begin{exercise}}
\ncom{\eexr}{\end{exercise}}  
\ncom{\bd}{\begin{definition}}
\ncom{\ed}{\end{definition}}  
\ncom{\brm}{\begin{remark}}   
\ncom{\erm}{\end{remark}}     
\ncom{\bal}{\begin{Algorithm}}
\ncom{\eal}{\end{Algorithm}}  
\ncom{\ol}{\overline}
\ncom{\wh}{\widehat} 

\ncom{\pf}{\noi {\bf Proof.  }}
\ncom{\eprf}{\noi {$\Box$}}
\ncom{\be}{\begin{enumerate}} 
\ncom{\ee}{\end{enumerate}}   
\ncom{\seq}{\subseteq}

\ncom{\zr}{\bf\textcolor{red}}
\ncom{\zb}{\bf\textcolor{blue}}
\ncom{\zg}{\bf\textcolor{green}}
\ncom{\zm}{\bf\textcolor{magenta}}

\definecolor{gold}{rgb}{0.85,.66,0}
\definecolor{gb}{rgb}{0, .5,.5}
\definecolor{rb}{rgb}{0.5, 0,.5}
\definecolor{Pink}{rgb}{1,0.75, 0.8}
\ncom{\zgb}{\bf\textcolor{gb}}
\ncom{\zrb}{\bf\textcolor{rb}}

\makeatletter
\newcommand{\single}{\let\CS=\@currsize\renewcommand{\baselinestretch}{1.5}\tiny\CS}
\newcommand{\oneandahalfspacing}{\let\CS=\@currsize\renewcommand{\baselinestretch}{1.5}\tiny\CS}
\newcommand{\doublespacing}{\let\CS=\@currsize\renewcommand{\baselinestretch}{1.6}\tiny\CS}
\newcommand{\double}{\let\CS=\@currsize\renewcommand{\baselinestretch}{3}\tiny\CS}

\def\notarro{{{\hbox{{\hspace*{-.02in}}$\rightarrow$ } }
{\hbox{$\!\!\!\!\!\!\!$}}{\raise 0.2ex 
\hbox{$\scriptscriptstyle{/}$}}}\hspace{.06in}}
\def \*{^{\mbox{$*$}}}

\newcommand\bnota{\begin{nota} }                         
\newcommand\enota{\end{nota} }                           

\newcommand{\beano}{\begin{eqnarray*}}
\newcommand{\eeano}{\end{eqnarray*}}

\usepackage{euscript}

\title{{\bf The perfect matching association scheme}}

\author{{\textcolor{black} {\bf Murali K. Srinivasan}} \\
{\em  {Department of Mathematics}}\\
{\em  {Indian Institute of Technology, Bombay}}\\
{\em  {Powai, Mumbai 400076, INDIA}}\\
{\bf  \texttt{mks@math.iitb.ac.in}}\\
{\bf  \texttt{murali.k.srinivasan@gmail.com}}
{\small Gelfand-Tsetlin vectors}
}

\begin{document}
\date{}
\maketitle
\begin{center}{\small{\em To the memory of Angel, Gandhi, and
Shadow}}\end{center}

\begin{abstract} 

We revisit the Bose-Mesner algebra of the perfect matching association
scheme. Our main results are: 
\begin{itemize}
\item An inductive algorithm, based on solving linear
equations, to compute the eigenvalues of the
orbital basis elements given the 
central characters of the symmetric groups.

\item Universal formulas,
as content evaluations of symmetric functions, for the eigenvalues of fixed
orbitals. 

\item An inductive construction of an eigenvector (the so called
first Gelfand-Tsetlin vector) in each eigenspace leading to a different
inductive algorithm (not using central characters) 
for the eigenvalues of the orbital basis elements.
\end{itemize}
\end{abstract}

{{\bf  \section {  Introduction }}}

In this paper we revisit the Bose-Mesner algebra of the perfect matching
association scheme.
The symmetric group $S_{2n}$ has a natural substitution action on the set
$\cM_{2n}$ of all perfect matchings in the complete graph
$K_{2n}$. The corresponding permutation representation of $S_{2n}$ on
$\C[\cM_{2n}]$ (the complex vector space with $\cM_{2n}$ as basis) 
is multiplicity free
and the (commutative) algebra  
$\cB_{2n} = \End_{S_{2n}}(\C[\cM_{2n}])$
is called the Bose-Mesner algebra of the perfect matching
association scheme. The eigenspaces of $\cB_{2n}$, in its left action on
$\C[\cM_{2n}]$, are indexed by even Young diagrams with $2n$ boxes 
(i.e., Young diagrams with $2n$ boxes
having an even number of boxes in every row)  
and the orbital basis 
of $\cB_{2n}$ is indexed by even partitions of $2n$ (i.e., partitions of
$2n$ with all parts even). The present work is motivated by the following
two results.

Diaconis and Holmes {\bf\cite{dh}} determined all the eigenvalues
of the orbital basis element of $\cB_{2n}$ indexed by the even partition
$(4,2^{n-2})$ of $2n$ (here $(4,2^{n-2})$
denotes the even partition with one part equal to 4 and $n-2$ parts
equal to $2$). 
We generalize this result to all fixed orbitals in Theorem
\ref{mks} below.

Godsil and Meagher {\bf\cite{gm,gm1}} and Lindzey {\bf\cite{l}} 
write down an eigenvector (using a
quotient argument) belonging to the eigenspace indexed by the even Young
diagram $(2n-2,2)$ with $2n$ boxes, yielding the eigenvalues of all orbital basis
elements on this eigenspace. 
We generalize this result by giving an inductive procedure to write down an
eigenvector in every eigenspace in Theorem \ref{mt1} below. This yields a
practical algorithm to compute the eigenvalues 
that we have implemented in Maple, see {\bf\cite{sr}}. The program computes, reasonably
efficiently, any given eigenvalue upto $\cB_{40}$. We were able to determine
the entire spectrum of the perfect matching derangement matrix in
$\cB_{2n}$,
upto $2n=40$ (see Problem 16.10.1 in {\bf\cite{gm}}).  

The rest of the introduction gives a more
detailed, although still informal, description of our results.

A partition (or a Young diagram) $\lambda$ is called {\em even} if all parts (or
all row lengths) of $\lambda$ are even. Clearly, $\lambda=(\lambda_1,\ldots
,\lambda_k) \mapsto 2\lambda=(2\lambda_1,\ldots ,2\lambda_k)$ is a bijection
between the set of all partitions of $n$ (or Young diagrams with $n$ boxes) 
and the set of all even partitions of $2n$ (or even
Young diagrams with $2n$ boxes). 
Let $\cP$ denote the set of all partitions
and $\cY$ denote the set of all Young diagrams (there is a unique partition
of $0$ and there is a unique Young diagram with $0$ boxes, both denoted
$(0)$). 
Let $\cP_n$ denote the set of all partitions of $n$
and let $\cY_n$ denote the set of all Young diagrams with $n$ boxes.
If $\lambda$ is
a partition of $n$ or if $\lambda$ is a Young diagram with $n$ boxes
we write $\lambda\vdash n$ and $|\lambda|=n$ (it will be clear from the
context whether a partition or a Young diagram is meant). 

Given a Young
diagram $\lambda$ with $n$ boxes, denote the (complex) irreducible
representation of $S_n$ parametrized by $\lambda$ by $V^{\lambda}$ and
denote the character of $V^{\lambda}$ by $\chi^{\lambda}$. 
For $\mu\vdash n$, denote the conjugacy class 
of permutations in $S_n$ of cycle type $\mu$
by $C_{\mu}$ and set $\chi^{\lambda}_{\mu}=\chi^{\lambda}(\pi)$, for (any)
$\pi\in C_{\mu}$. We let $k_{\mu}\in \C[S_n]$ (= the group algebra of $S_n$) 
denote the sum of elements in $C_{\mu}$.

Let $Z[\C[S_n]]$ denote the center of the group algebra of $S_n$.
Then $Z[\C[S_n]]$ is a semisimple commutative
algebra of dimension $p(n)$, the number of
partitions of $n$, with $\{k_\mu\;|\;\mu\vdash n\}$ as a basis.  
The eigenspaces of this algebra, in its left action on $\C[S_n]$, 
are the isotypical
components of $V^\lambda,\;\lambda\vdash n$ in $\C[S_n]$. Let
$\hat{\phi}^{\lambda}_\mu$ denote the eigenvalue of $k_\mu$ on the isotypical
component of $V^\lambda$.
By taking traces we see that    
\beq \label{equ}
\hat{\phi}^{\lambda}_{\mu} &=&
\frac{|C_{\mu}|\chi^{\lambda}_{\mu}}{\dim(V^{\lambda})}.
\eeq
We call $\hat{\phi}^{\lambda}_{\mu}$ a {\em central character}. 
As there are well known explicit formulas for $|C_\mu|$ and
$\dim(V^\lambda)$ we may regard  
$\hat{\phi}^\lambda_\mu$ and $\chi^\lambda_\mu$ as being equivalent from the
point of view of computing them.
It can be easily 
shown that $\hat{\phi}^\lambda_\mu$ is an integer.

We now define an analog of $Z[\C[S_n]]$.
We have the following basic result (see  
{\bf\cite{bu,jk,m,s,st}}): there is a $S_{2n}$-linear isomorphism
\beq \label{fr}
\C[\cM_{2n}] & \cong & \oplus_{\lambda\vdash n} V^{2\lambda}.
\eeq

Let $\cB_{2n} = \End_{S_{2n}}(\C[\cM_{2n}])$. Since $\C[\cM_{2n}]$ is
multiplicity free, $\cB_{2n}$ is a semisimple commutative algebra called the {\em
Bose-Mesner algebra of the perfect matching association scheme}. Its
dimension is also $p(n)$. 

From (\ref{fr}) above we have that the common eigenspaces of $\cB_{2n}$, in its
left action on $\C[\cM_{2n}]$, are ($S_{2n}$-isomorphic to)
$V^{2\lambda},\;\lambda\vdash n$.
The orbits of the diagonal action of
$S_{2n}$ on $\cM_{2n}\times \cM_{2n}$, and thus the orbital basis of
$\cB_{2n}$, can be shown to be indexed by even partitions of $2n$ (see
Section 2). Given $\mu\vdash n$, let $N_{2\mu}$ denote the orbital basis
element of $\cB_{2n}$ indexed by the even partition $2\mu$ and let
$\hat{\theta}_{2\mu}^{2\lambda},\;\lambda,\mu\vdash n$, 
denote the eigenvalue (which can be shown to be an integer, see Section 2)
of $N_{2\mu}$ on 
$V^{2\lambda}$. We refer to the $\hat{\theta}_{2\mu}^{2\lambda}$ as the {\em
eigenvalues} of $\cB_{2n}$. We think of $\hat{\theta}^{2\lambda}_{2\mu}$ as an analog of
$\hat{\phi}^{\lambda}_{\mu}$.

In Section 3 we address the following question:
assuming the central characters of $S_n$ as given, how can we
calculate the eigenvalues of the Bose-Mesner algebra.  
We give a recursive combinatorial algorithm for this task:
we show that we can inductively compute the eigenvalues
of $\cB_2,\cB_4,\ldots ,\cB_{2n}$ from the central characters of
$S_2,S_4,\ldots ,S_{2n}$ by solving systems of linear equations.

Let $\hat{\Theta}(2n)$ denote the eigenvalue table of $\cB_{2n}$, i.e.,
$\hat{\Theta}(2n)$ is the $\cY_n\times \cP_n$ matrix with entry in row
$\lambda$, column $\mu$ given by $\hat{\theta}^{2\lambda}_{2\mu}$.

\bt \label{mt} Assume given the central characters of
$S_2, S_4, \ldots , S_{2n}$ and the eigenvalues of $\cB_2, \cB_4, \ldots
,\cB_{2n-2}$. There is an algorithm that determines the eigenvalues of
$\cB_{2n}$ by solving nonsingular systems of linear equations with
coefficient matrices  
$\hat{\Theta}(2), \hat{\Theta}(4), \ldots , \hat{\Theta}(2n-2)$ and with
right hand sides determined by the central characters of $S_4, S_6, \ldots
,S_{2n}$. 

Thus we can inductively compute the eigenvalues of the Bose-Mesner algebra
from the central characters of the symmetric groups 
by solving linear equations.
\et

Theorem \ref{mt}, when combined with the work of Corteel, Goupil, and
Schaeffer {\bf\cite{cgs}} and Garsia {\bf\cite{g}} expressing central
characters (at fixed conjugacy classes) as content evaluations of symmetric
functions, yields similar formulas for the eigenvalues of fixed orbital basis
elements. Let us explain this. First we introduce notation concerning fixed
classes and symmetric functions.

Let $\cP(2)$ denote the set of  
partitions with all parts $\geq 2$. Note that
the unique partition of $0$ belongs to $\cP(2)$. 
For $\mu\in \cP(2)$, let
$\ol{\mu}$ be the partition of $|\mu|-\ell(\mu)$ ($\ell(\mu)$ = number of
parts of $\mu$) obtained by subtracting $1$
from every part of $\mu$. The map $\cP(2) \rar \cP$ given by $\mu \mapsto
\ol{\mu}$ is clearly a bijection. Let $\cP(2,n)$ denote the set of all
$\mu\in\cP(2)$ with $|\mu|\leq n$.

By a nontrivial cycle of a permutation we mean a cycle of length $\geq 2$.
Given $\mu\in \cP(2)$ and $n\geq 1$, define $c_{\mu}(n)$ to be 
element of 
$Z[\C[S_n]]$ given by the sum of all permutations $\pi$ in $S_n$ 
that have $\mu$ as the
partition determined by the lengths of the nontrivial cycles of $\pi$. Thus, 
$c_\mu(n)$ is $0$ if $n<|\mu|$ and is equal to 
$k_{(\mu,1^{n-|\mu|})}$ if $n\geq |\mu|$ (here $(\mu,1^{n-|\mu|})$
denotes the partition of $n$ obtained by adding, to $\mu$, $n-|\mu|$ parts
equal to 1). 
In this notation, $c_{(3)}(n)$ denotes the conjugacy class sum of 3-cycles
in $\C[S_n]$ (which is automatically zero if $n=1,2$), $c_{(0)}(n)$
denotes the identity element of $\C[S_n]$, and 
$\{c_\mu(n)\;|\;\mu\in\cP(2,n)\}$ is a basis of $Z[\C[S_n]]$.

Given $\mu\in \cP(2)$ and $\lambda\in \cY$, define $\phi^\lambda_\mu$ to be
the eigenvalue of $c_\mu(|\lambda|)$ on $V^\lambda$. That is, if
$\lambda$ has $n$ boxes, $\phi^\lambda_\mu$ is equal to
$\hat{\phi}^\lambda_{(\mu,1^{n-|\mu|})}$ if $n\geq |\mu|$ and is equal to
$0$ if $n<|\mu|$.

Similarly, given $\mu\in \cP(2)$ and $n\geq 1$, define $M_{2\mu}(2n)$ to be the element of
$\cB_{2n}$ given as follows: it is equal to  the orbital basis element
$N_{2(\mu,1^{n-|\mu|})}$ if $n\geq |\mu|$ and it is $0$ if $n<|\mu|$.
For instance, if $\mu = (3,2,1,1)\vdash 7$ and $\tau =(3,2)$ we can write
the element $N_{2\mu}$ of $\cB_{14}$ as $M_{2\tau}(14)$. The orbital basis
of $\cB_{2n}$ can be written as $\{M_{2\tau}(2n)\;|\;\tau\in\cP(2,n)\}$.

Given $\mu\in \cP(2)$ and $\lambda\in \cY$, define $\theta^{2\lambda}_{2\mu}$ to be
the eigenvalue of $M_{2\mu}(2|\lambda|)$ on $V^{2\lambda}$. That is, if
$\lambda$ has $n$ boxes, $\theta^{2\lambda}_{2\mu}$ is equal to
$\hat{\theta}^{2\lambda}_{2(\mu,1^{n-|\mu|})}$ if $n\geq |\mu|$ and is equal to
$0$ if $n<|\mu|$.

We think of $\hat{\phi}^\lambda_\mu$ 
and $\hat{\theta}^{2\lambda}_{2\mu}$ as
functions of $\lambda, \mu\vdash n$, for fixed $n$ . While considering
$\phi^{\lambda}_{\mu}$ and $\theta^{2\lambda}_{2\mu}$, we regard $\mu$ as
fixed, and think of $\phi^\lambda_\mu$, $\theta^{2\lambda}_{2\mu}$ 
as functions on $\cY$.

The {\em content}
$c(b)$ of a box $b$ of a Young diagram $\lambda$ is its
$y$-coordinate minus its $x$-coordinate (our convention for drawing Young
diagrams is akin to writing down matrices with $x$-axis running downwards
and $y$ axis running to the right). Thus the content of the boxes in the
first row (from left to right) are $0,1,2,\ldots $, in the second row are
$-1,0,1, \ldots$, and so on. We denote by 
$c(\lambda)$ the multiset of contents 
of all the boxes of $\lambda$. 
So $c(\lambda)$ has (multiset) cardinality $|\lambda|$.

Let $\Lambda[t]$ denote the algebra, over $\Q[t]$, of symmetric functions in
$\{x_1,x_2,x_3,\ldots \}$. Define $p_0=1$ and
$p_n=\sum_{i}x_i^n,\;n\geq 1$. For $\lambda\in\cP$ the {\em power sum
symmetric function} $p_{\lambda}$ is defined as follows:
\beqn
p_{\lambda} &=& p_{\lambda_1}p_{\lambda_2}\cdots\;\;\;\;\;\mbox{if
}\lambda=(\lambda_1,\lambda_2,\ldots ).
\eeqn 
The set $\{p_{\lambda}\;|\;\lambda\in\cP\}$ is a $\Q[t]$-module basis of
$\Lambda[t]$ ({\bf\cite{cst2,m,p,sa,st}}).

Given $f\in \Lambda[t]$ and $\lambda\in\cY$ with $n$ boxes we 
define the {\em content evaluation} $f(c(\lambda))$ to be the 
rational number obtained from $f$ by
setting $t=n,\;x_i=0\mbox{ for }i>n$, and 
$$\{x_1,x_2,\ldots ,x_n\}=\mbox{ (the multiset) }c(\lambda).$$
Note that this definition makes sense as $f$ is symmetric.

Frobenius proved that the central character at the conjugacy class of
transpositions is given by content evaluation of the symmetric function
$p_1\in \Lambda[t]$ and Ingram proved that the central character at the
conjugacy class of 3-cycles is given by content evaluation of the
symmetric function $p_2-\frac{t(t-1)}{2} \in \Lambda[t]$ (see {\bf\cite{cgs}}).
These are {\em universal formulas} (i.e., independent of 
$\lambda$) made precise as
follows:
\beqn
\phi^{\lambda}_{(2)} &=& p_1(c(\lambda)) = \mbox{ 
Sum of contents of all boxes of $\lambda$ },\;\;\lambda\in \cY,\\
\phi^{\lambda}_{(3)} 
&=& \left(p_2-\frac{t(t-1)}{2}\right)(c(\lambda))\\
&=&\mbox{ Sum of squares of contents of all boxes of $\lambda$ } 
- \frac{|\lambda|(|\lambda|-1)}{2},\;\;\lambda\in \cY.
\eeqn
Note that $\phi^\lambda_{(3)}$ is 0 when
$|\lambda|=1,2$. 
These formulas can be generalized to all fixed conjugacy classes.

For each $\mu\in\cP(2)$, it is shown in {\bf\cite{cgs}} that  
there is a symmetric
function $W_{\mu}\in \Lambda[t]$ such that
$\{W_{\mu}\;|\;\mu\in\cP(2)\}$ is a $\Q[t]$-module basis of $\Lambda[t]$ and,
for all $\mu\in \cP(2),\;\lambda\in \cY$,
\beqn
\phi^{\lambda}_{\mu} &=& W_\mu(c(\lambda)).
\eeqn
An algorithm to compute $W_\mu$ is given in {\bf\cite{g}}.
We motivate and discuss this result in Section 4. 

Diaconis and Holmes {\bf\cite{dh}} observed, using Frobenius' result, 
that the eigenvalues of the
orbital basis element of $\cB_{2n}$
corresponding to 4-cycles (i.e., the even partition $(4,2^{n-2})$) are given by content
evaluation of the symmetric function $\frac{p_1}{2}- \frac{t}{4}\in
\Lambda[t]$, i.e.,
\beqn
\theta^{2\lambda}_{2(2)} &=& \left(\frac{p_1}{2}-
\frac{t}{4}\right)(c(2\lambda))\\\\ 
 &=& \frac{\mbox{ 
Sum of contents of all boxes 
of $2\lambda$ }}{2} - \frac{2|\lambda|}{4},\;\;\;\;\;\;\lambda\in \cY.
\eeqn
Note that $\theta^{2\lambda}_{2(2)}=0$ when $|\lambda|=1$.
This can be generalized to all fixed orbital basis elements.

In Section 4, we show that the algorithm of Theorem \ref{mt} 
converts the basis $\{W_\mu\}$ of $\Lambda[t]$ into another basis $\{E_\mu\}$ of
$\Lambda[t]$ with the following property.

\bt \label{mks}
For each $\mu\in\cP(2)$ there is a symmetric
function $E_{\mu}\in \Lambda[t]$ such that

\noi (i) $\{E_{\mu}\;|\;\mu\in\cP(2)\}$ is a $\Q[t]$-module basis of
$\Lambda[t]$.

\noi (ii) For all $\mu\in\cP(2)$ and $\lambda\in\cY$, we have
\beqn
\theta^{2\lambda}_{2\mu} &=& E_{\mu}(c(2\lambda)).
\eeqn
\et
Information about the coefficients in the expansion of $W_\mu$ and
$E_\mu$ in the power sum basis is given in Section 4. Example \ref{wep} in
Section 4 lists these symmetric functions for $|\mu|\leq 4$.

One method for computing the eigenvalues $\{\theta_i\}$ of a real symmetric
matrix $N$ is to write down eigenvectors $\{v_i\}$, one in each eigenspace,
and then to solve for $\theta_i$ in the equation $Nv_i=\theta_iv_i$. 
In Section 5 we use this method to give a different  
inductive algorithm (not using the characters or central characters of
$S_n$) for computing the
eigenvalues $\hat{\theta}^{2\lambda}_{2\mu}$ of $\cB_{2n}$. 

Every $S_n$-irreducible $V^\lambda$ has a canonically
defined basis, determined upto scalars, and called the Gelfand-Tsetlin (GZ) basis. We
systematically choose one of these basis vectors  and call it the first GZ
vector (see Sections 4 and 5 for definitions). Denote the first GZ vector of
$V^\lambda$ by $v_\lambda$. 
Let $\lambda'\in \cY_{n+1}$ with $\lambda =
\lambda' - \{\mbox{ last box in the last row of $\lambda'$ }\}$. Let
$v_{2\lambda}$ denote the first GZ vector of the eigenspace $V^{2\lambda}$
of $\cB_{2n}$. Then there is a simple expression for $v_{2\lambda'}$ in terms
of $v_{2\lambda}$ (see Section 5). 
The simplest nontrivial case of this occurs when
$\lambda'=(n,1)$. Here $\lambda=(n)$ and $V^{2(n)}$ is the trivial
representation giving $v_{2\lambda}=\sum_{A\in \cM_{2n}} A$. In this case
the eigenvector vector $v_{2\lambda'}$ coincides with that written down by
Godsil and Meagher {\bf\cite{gm,gm1}} and Lindzey {\bf\cite{l}} (using a
quotient argument).

Of course, explicitly writing down these
vectors is inefficient since $v_{2\lambda}$ lives in a space
of dimension $(2n-1)!!=1\cdot 3\cdot 5 \cdots \cdot (2n-1)$. 
However, we use this expression implicitly to
give an algorithm that works with only the rows of $\hat{\Theta}(2n)$.
Note that a row of $\hat{\Theta}(2n)$ has only $p(n)$ components, which is
subexponential and is only moderately large for small values of $n$ 
(for example, compare $p(13)=101$ with $25!!=7905853580625$).
\bt \label{mt1}
Let $\lambda'\in \cY_{n+1}$ with $\lambda =
\lambda' - \{\mbox{ last box in the last row of $\lambda'$ }\}$. 
Assume that
the row of $\hat{\Theta}(2n)$ indexed by $\lambda$, i.e., the vector
$(\hat{\theta}^{2\lambda}_{2\mu})_{\mu\vdash n}$, is known.

There is an algorithm to determine
$(\hat{\theta}^{2\lambda'}_{2\mu'})_{\mu'\vdash n+1}$, i.e.,
the row of $\hat{\Theta}(2n+2)$ indexed by $\lambda'$.

\et
The eigenvector approach also applies to the central characters and in
Section 5 we give a virtually identical
inductive algorithm (not using irreducible characters) to compute  
$\hat{\phi}^\lambda_\mu$. Although this method of computing the central
characters is not as efficient as the one based on (\ref{equ}) (since the
irreducible characters can be very efficiently calculated), 
it further brings out the
essential analogy between $\hat{\theta}^{2\lambda}_{2\mu}$ 
and $\hat{\phi}^\lambda_\mu$. A simple recursive implementation of these
algorithms in Maple is given in {\bf\cite{sr}}. This program calculates,
reasonably quickly,
$\hat{\phi}^\lambda_\mu$ and $\hat{\theta}^{2\lambda}_{2\mu}$, for
$|\lambda|=|\mu| \leq 20$.

Finally, we would like to add a terminological remark. The Bose-Mesner
algebra $\cB_{2n}$ is isomorphic to the Hecke algebra (also called the
double coset algebra) of the Gelfand pair
$(S_{2n},H_n)$, where $H_n$ is the hyperoctahedral group (see Example 5 
of Chapter VII.2 of {\bf\cite{m}}), and the two settings are equivalent. In
this paper we adopt the perfect matching point of view.

{{\bf  \section {The $S_n$-module $\C[\cM_n]$}}}

The regular modules $\C[S_n]$ have the following recursive
structure
\beq \label{rrs}
\ind_{\;S_n}^{\;S_{n+1}}(\C[S_n])\cong  \C[S_{n+1}].
\eeq
The modules $\C[\cM_{2n}]$ have a similar recursive structure. Informally,
we can say that
the induction happens at every other step and we do nothing in between (see
items (v) and (vi) of Lemma \ref{mfree} below). This 
is best brought out by
simultaneously considering the odd case, i.e., the action of $S_{2n+1}$ on 
near perfect matchings (= matchings with $n$ edges) of $K_{2n+1}$. This idea
is implicit in the detailed proof of (\ref{fr}) given in Chapter 43 of
Bump's book {\bf\cite{bu}} (also see {\bf\cite{jk,s}}) but it is
useful to make it explicit as it simplifies certain technicalities and
also suggests an approach to writing down the
eigenvectors of $\cB_{2n}$ in
Section 5. We adopt a uniform notation for both the even and odd cases.

Let $\cbP_n$ denote the set of all even partitions of
$n$, if $n$ is even, 
or the set of all near even partitions of $n$ (i.e., exactly 
one part odd), if $n$ is odd. 
Let $\cbY_n$ denote the set 
of all even Young diagrams with
$n$ boxes, if $n$ is even, 
or the set of all near even Young diagrams with $n$ boxes (i.e., exactly 
one row length odd), if $n$ is odd.

Let $\cM_n$ denote the set of all maximum matchings in $K_n$ (i.e., perfect
matchings if $n$ is even and near perfect matchings if $n$ is odd).
Given $A,B \in \cM_n$, let $d(A,B)$ be the partition whose parts
are the number of vertices in the spanning subgraph of $K_n$ with edge set
$A\cup B$. It is easily seen that $d(A,B)\in \cbP_n$.

For $\mu \in \cbP_n,\;A\in \cM_n$ define
\beqn
\cM(A,\mu) &=& \{ B\in \cM_n\;|\; d(A,B)=\mu\},
\eeqn
and define a linear operator
\beqn
N_\mu: \C[\cM_n]\rar \C[\cM_n]
\eeqn
by setting, for $A\in \cM_n$,
$$N_\mu(A)=\sum_{B\in \cM(A,\mu)}\;B.$$

The symmetric group $S_n$ has a natural action on $\cM_n$ and this gives
rise to the $S_n$-module $\C[\cM_n]$. We have the diagonal action of $S_n$
on $\cM_n\times \cM_n$. Set $\cB_n = \End_{S_n}(\C[\cM_n])$.

For $n$ odd, given $A\in \cM_n$ we denote by $v(A)$ the unique vertex of
$K_n$ that is not the endpoint of any edge in $A$. An edge connecting
vertices $i$ and $j$ will be denoted $[i,j]$ (or $[j,i]$). The following result
collects together basic properties of the $S_n$-action on $\cM_n$.

\bl \label{mfree} Let $n$ be a positive integer.

\noi (i) $(A,B),\;(C,D)\in \cM_n\times \cM_n$ are in the same $S_n$-orbit if
and only if $d(A,B)=d(C,D)$.

\noi (ii) The set $\{ N_\mu\;|\; \mu\in \cbP_n\}$ is a basis of
$\cB_n$.

\noi (iii) $(A,B),\,(B,A)$ are in the same $S_n$-orbit, for all $(A,B)\in
\cM_n\times \cM_n$.

\noi (iv) The $S_n$-module $\C[\cM_n]$ is multiplicity free.

\noi (v) Assume $n$ is odd. We have an $S_n$-module isomorphism
(treating $S_n$ as the subgroup of $S_{n+1}$ fixing $n+1$) 
$$\C[\cM_n] \cong \res_{\;S_n}^{\;S_{n+1}}(\C[\cM_{n+1}])$$
given by $A\mapsto A \cup \{[v(A),n+1]\},\;A\in \cM_n$. 

\noi (vi) Assume $n$ is even. We have an $S_{n+1}$-module isomorphism
$$\ind_{\;S_n}^{\;S_{n+1}}(\C[\cM_n])\cong  \C[\cM_{n+1}].$$
\el  

\pf (i) This is clear.

\noi (ii) 
This follows from part (i) by a standard result (see {\bf\cite{cst2,
gm}}). This basis is called the {\em orbital basis} of the commutant.

\noi (iii) Follows from part (i).

\noi (iv) This follows from part (iii) by a standard result (see
{\bf\cite{cst2,gm}}).

\noi (v) This is clear.

\noi (vi) Consider the disjoint union given by coset decomposition 
\beqn S_{n+1}&=&S_n \cup (1\; n+1)S_n \cup \cdots \cup (n\;n+1)S_n.
\eeqn
We think of $\ind_{\;S_n}^{\;S_{n+1}}(\C[\cM_n])$ as the (left)
$\C[S_{n+1}]$-module $\C[S_{n+1}] \otimes_{\C[S_n]} \C[\cM_n]$ with
basis $\{(i\;n+1)\otimes A\;:\;1\leq i \leq n+1,\;A\in \cM_n\}$ (here
$(n+1\;n+1) = 1$, the identity element of $S_{n+1}$).

Define a bijective linear map
$f:\ind_{\;S_n}^{\;S_{n+1}}(\C[\cM_n])\rar \C[\cM_{n+1}]$ by
\beqn
f((i\;n+1)\otimes A) &=& (i\;n+1)\cdot A,\;\;\;1\leq i \leq n+1,\;A\in \cM_n.
\eeqn
Fix $1\leq i \leq n+1$ and $A\in \cM_n$. Let $\tau \in S_{n+1}$. Set
$j=\tau(i)$ and write
$\tau (i\;n+1)=(j\;n+1)\tau'$ where $\tau'=(j\;n+1)\tau (i\;n+1)$. Note that
$\tau'(n+1)=(n+1)$. Then
\beqn
f(\tau\cdot((i\;n+1)\otimes A))&=&f((j\;n+1)\otimes (j\;n+1)\tau(i\;n+1)\cdot
A) \\
&=& (j\;n+1)\cdot ((j\;n+1)\tau(i\;n+1)\cdot A)\\
&=&\tau\cdot f((i\;n+1)\otimes A).
\eeqn
Thus, $f$ is an $S_{n+1}$-module isomorphism. \eprf

Parts (ii) and (iv) of Lemma \ref{mfree} show that the eigenvalues of
$N_\mu$ are integers using the following standard argument (and the
fact that the irreducible characters of $S_n$ are integer valued).

\bl \label{inteig}
Let a finite group $G$ act on a finite set $X$ and for, $g\in G$, let
$\rho(g)$ denote the $X\times X$ permutation matrix corresponding to the
action of $g$ on $X$. Let $A$ be a $X\times X$
matrix with integer entries that commutes with the action of $G$ on $X$,
i.e., $A\rho(g)=\rho(g)A$ for all $g\in G$. 
Assume that 

\noi (i) The permutation representation of $G$ on $\C[X]$ is multiplicity free.  

\noi (ii) The character of every $G$-irreducible appearing in $\C[X]$ is
integer valued.

\noi Then the eigenvalues of $A$ are integral.
\el
\pf Write 
$$\C[X] = V_1 \oplus \cdots \oplus V_t,$$
where $V_1,\ldots ,V_t$ are nonisomorphic irreducible $G$-submodules of
$\C[X]$. Let $\chi_i$ be the character of $V_i$.

Let $\lambda$ be an eigenvalue of $A$. 
By Schur's lemma, every $V_j$ is contained in an eigenspace of $A$. 
Thus the eigenspace of $\lambda$ is a direct sum of some of the $V_j$'s. 
Say $V_i$ is contained in the eigenspace of $\lambda$.

The
$G$-linear projection $\C[X] \rar \C[X]$ onto $V_i$ is given by
$$ v\mapsto \frac{\dim V_i}{|G|} \sum_{g\in G} \ol{\chi_i(g)}\;g\cdot v$$
Since $\chi_i$ is integer valued the matrix of the projection above (in the
standard basis $X$) has rational entries and thus there is an eigenvector
for $\lambda$ with rational entries. Since $A$ is integral it follows that
$\lambda$ is a rational number and since it is also an algebraic integer
(being an eigenvalue of an integer matrix) it follows that $\lambda$ is an
integer. \eprf

The recursive structure of the modules $\C[\cM_n]$ given by parts (v) and
(vi) of Lemma \ref{mfree}, together with the branching rule, yields
a proof of (\ref{fr}). This part of the proof, which we include for
completeness, is essentially the same as in {\bf\cite{bu}}. 
Let us first recall the branching rule.

A fundamental
result (see {\bf\cite{cst2,jk,ov, p, sa}}) in the representation theory of the
symmetric groups states that there is a unique assignment,
denoted $\lambda
\mapsto V^{\lambda}$, which associates to each Young diagram $\lambda$ an
equivalence class $V^{\lambda}$ of irreducible $S_{|\lambda|}$-modules
(we also let $V^{\lambda}$ 
denote an irreducible $S_n$-module in the corresponding equivalence
class)
such
that properties (a) and (b) below are satisfied:

(a) {\em Initialization:}  $V^{(2)}$ is the trivial
representation of $S_2$ 
and $V^{(1,1)}$ is the sign representation of $S_2$
(here $(2)$, respectively $(1,1)$, denotes the Young diagram with a single
row of two boxes, respectively a single column of two boxes).

(b) {\em Branching rule:}
Given $\mu\in \cY$, we denote by
$\mu^-$ the set of all Young diagrams obtained from $\mu$ by
removing a box corresponding to 
one of the inner corners in the Young diagram $\mu$.
For $n\geq 2$, given $\lambda\in \cY_n$, consider the irreducible
$S_n$-module $V^{\lambda}$.
Viewing $S_{n-1}$ as the subgroup
of $S_n$ fixing $n$
we have an $S_{n-1}$-module isomorphism
\beq \label{brres}
&\res^{\;S_n}_{\;S_{n-1}}(V^{\lambda})\cong  \bigoplus_{\mu\in\lambda^-}
V^{\mu}.&
\eeq

It is a consequence of properties (a) and (b) above that
$\{V^{\lambda}\;|\;\lambda\in \cY_n\}$ is a complete set of pairwise
inequivalent irreducible representations of $S_n$. Another consequence is
that, for any $n$, the Young diagram consisting of a single row of $n$ boxes
(respectively, a single column of $n$ boxes)    
corresponds to the trivial representation of $S_n$ (respectively, the sign
representation of $S_n$). 

Given $\mu\in \cY$, we denote by
$\mu^+$ the set of all Young diagrams obtained from $\mu$ by
adding a box corresponding to one of the outer 
corners in the Young diagram $\mu$.
For $n\geq 1$, given $\lambda\in \cY_n$, consider 
the irreducible $S_{n}$-module $V^{\lambda}$. 
By Frobenius reciprocity, the branching rule
can be equivalently stated as
\beq \label{brind}
&\ind_{\;S_n}^{\;S_{n+1}}(V^{\lambda})\cong  \bigoplus_{\mu\in\lambda^+}
V^{\mu}.&
\eeq

\bt \label{ep} Let $n$ be a positive integer. There is a $S_n$-linear
isomorphism
\beqn \C[\cM_n] \cong \bigoplus_{\lambda\in \cbY_n} V^{\lambda}.
\eeqn
\et 
\pf The proof is by induction on $n$, the cases $n=1,2$ being clear. Let
$n\geq 3$ and consider the following two cases.

\noi (i) $n$ is odd: This easily follows from the induction hypothesis, 
Lemma \ref{mfree} (iv), (vi), and the branching rule.

\noi (ii) $n$ is even: Let $V^{\lambda}$, $\lambda\in \cY_n$ occur in
$\C[\cM_n]$ and assume that $\ell(\lambda)\geq 3$. Suppose that all rows of
$\lambda$ are not of even length. Then, since $n$ is even, we can find an
inner corner of $\lambda$ such that deleting the corresponding box leaves a
Young diagram with at least two rows of odd length. By Lemma \ref{mfree} (v)
and the branching rule,
this contradicts the induction hypothesis (for $n-1$). Thus,
$V^{\lambda}$ cannot occur in $\C[\cM_n]$.

Define Young diagrams $\lambda_k=(n-k,k),\;0\leq k \leq n/2$. Note that
$\lambda_0,\ldots ,\lambda_{n/2}$ are all the Young diagrams with at most
two rows. We shall show, by induction on $k$, that $V^{\lambda_k},\;0\leq k\leq
n/2$ occurs in $\C[\cM_n]$ if and only if $k$ is even. Now $V^{\lambda_0}$ is
the trivial representation and thus occurs in permutation representation
$\C[\cM_n]$. Assume, inductively, that our claim has been proven for
$V^{\lambda_0},\ldots ,V^{\lambda_{t-1}}$ and consider $V^{\lambda_t}$.
Suppose $t$ is even. By the main induction hypothesis on $n$,
$V^{(n-t,t-1)}$ occurs in $\C[\cM_{n-1}]$. By Lemma \ref{mfree} (v) and the
branching rule, one of $V^{(n-t,t-1,1)},\;V^{(n-t+1,t-1)},\;V^{(n-t,t)}$
must occur in $\C[\cM_n]$. The first cannot occur by the paragraph above,
the second cannot occur by the secondary induction hypothesis on $k$, and so
the third must occur. Now suppose that $t$ is odd and that $V^{(n-t,t)}$
occurs in $\C[\cM_n]$. Then, since $V^{(n-t+1,t-1)}$ occurs in $\C[\cM_n]$
(by the secondary  induction hypothesis on $k$), $V^{(n-t,t-1)}$ will occur
at least twice in $\C[\cM_{n-1}]$ contradicting its multiplicity freeness.
Thus the claim on $V^{\lambda_k},\;0\leq k \leq n/2$ is established.

What we have shown so far implies that if $V^{\lambda},\;\lambda\in \cY_n$
occurs in $\C[\cM_n]$ then all rows of $\lambda$ must have even length. Since,
by the branching rule,  
$\res^{\;S_n}_{\;S_{n-1}}(V^{\lambda})$ and
$\res^{\;S_n}_{\;S_{n-1}}(V^{\mu})$, for $\lambda,\;\mu\in \cbP_n,
\mu\not=\lambda$ can have no irreducibles in common, the result follows from
the induction hypothesis and Lemma \ref{mfree} (v). \eprf

{{\bf  \section {Eigenvalues and (class-coset) intersection numbers}}}

Assuming the
central characters of $S_2, S_4, \ldots , S_n$ as given, we
show in this section that we can compute
the eigenvalues of $\cB_{2n}$ by solving linear equations.

We begin by recalling, without proof, the following classical formula for the
eigenvalues of $\cB_{2n}$ that appears in Bannai and Ito {\bf\cite{bi}} 
(see page 179), Hanlon, Stanley, and Stembridge
{\bf\cite{hss}} (see equation (3.3) of Lemma 3.3)  and in
Godsil and Meagher {\bf\cite{gm}} (see Lemma 13.8.3). It is proved by
writing down the primitive idempotents of $\cB_{2n}$ and then 
expanding the orbital basis in terms of these. Another paper, using Jack
symmetric functions, on the
eigenvalues of $\cB_{2n}$ is Muzychuk {\bf\cite{muz}}.

Denote by $I$ the perfect matching $\{[1,n+1],[2,n+2],\ldots ,[n,2n]\}$
of $K_{2n}$. If $\mu$ is a
partition with $m_i$ parts equal to $i$ we set
$z_{\mu}=1^{m_1}m_1!\;2^{m_2}m_2!\;3^{m_3}m_3!\cdots$. 

\bt {\bf\cite{bi,hss,gm}} 
Let $\lambda, \mu \vdash n$. Fix $A\in \cM_{2n}$ with
$d(I,A)=2\mu$. 
Then
\beqn 
\hat{\theta}^{2\lambda}_{2\mu} &=&
\frac{1}{2^{\ell(\mu)}z_{\mu}}\left\{\sum_{\pi\in S_{2n},\;\pi\cdot I = A}
\chi^{2\lambda}(\pi)\right\}.\;\;\;\Box
\eeqn 
\et

The formula above has $2^nn!$ terms on the right hand side. We can group terms by cycle type to
reduce this number.

Let $\mu \vdash n$. Fix $A\in \cM_{2n}$ with $d(I,A)=2\mu$.
For $\tau\vdash 2n$, define 
$$m(\tau,2\mu) = |C_{\tau} \cap \{\pi\in S_{2n}\;|\; \pi\cdot I = A\}|,$$
i.e., $m(\tau, 2\mu)$ is the number of permutations in $S_{2n}$ 
of cycle type $\tau$ taking $I$
to $A$ (this number is clearly independent of $A$ as long as $d(I,A)=2\mu$). 
We refer to the $m(\tau,2\mu)$ as the {\em (class-coset) intersection
numbers} of $\cB_{2n}$ (being the cardinality of the 
intersection of a conjugacy class with a coset of the subgroup fixing $I$).

We thus have the following formula which has only $p(2n)$ terms
\beq \label{hssf}
\hat{\theta}^{2\lambda}_{2\mu} &=&
\frac{1}{2^{\ell(\mu)}z_{\mu}}\left\{\sum_{\tau\vdash
2n}m(\tau,2\mu)\chi^{2\lambda}_{\tau}\right\}.
\eeq 
There is, however, no simple formula for $m(\tau,2\mu)$. 
Thus, in the identity (\ref{hssf}) above, the characters of $S_{2n}$ 
are known but we have two sets of unknowns: eigenvalues of $\cB_{2n}$ 
and the intersection numbers of $\cB_{2n}$.
The idea 
of the present approach is the following bootstrap procedure: 
assume we have calculated the intersection numbers 
of $\cB_2,\ldots ,\cB_{2n-2}$. Then,

(i)  In Theorems \ref{ev1} and \ref{ev2} below we show that 
the eigenvalues of $\cB_{2n}$ can be found from the central characters of
$S_{2n}$ and the intersection numbers of $\cB_{2},\ldots ,\cB_{2n-2}$.

(ii) In Lemma \ref{mtrick} below we show that we  
can find the intersection numbers of $\cB_{2n}$ 
from the central characters of $S_{2n}$ and the eigenvalues of
$\cB_{2n}$ by solving linear equations.

For $\tau\vdash n, \mu\vdash 2n$ define column vectors of length $p(n)$
$$ \hat{\phi}_{\mu} = (\hat{\phi}^{2\lambda}_{\mu})_{\lambda\vdash n} 
\mbox{ and }
\hat{\theta}_{\tau} = (\hat{\theta}^{2\lambda}_{2\tau})_{\lambda\vdash n}.
$$
Note that $\hat{\theta}_\tau$ is the
column of $\hat{\Theta}(2n)$ indexed by $\tau$. We have
\bl \label{mtrick}
Let $\mu\vdash 2n$. Then
\beqn 
\hat{\phi}_{\mu} &=& \sum_{\tau\vdash n} m(\mu,2\tau)\hat{\theta}_{\tau},
\eeqn
i.e., defining the column vector $m(\mu)=(m(\mu,2\tau))_{\tau\vdash n}$ we have
\beqn
\hat{\phi}_\mu &=& \hat{\Theta}(2n)m(\mu).
\eeqn
\el
\pf Consider the element $k_{\mu}\in Z[\C[S_{2n}]]$. Then
\beqn 
k_\mu \cdot I &=& \sum_{\tau\vdash n} m(\mu,2\tau) N_{2\tau}(I). \eeqn
It follows that the actions of $k_\mu$ and 
$\sum_{\tau\vdash n} m(\mu,2\tau) N_{2\tau}$
on $\C[\cM_{2n}]$ are identical. The eigenvalue of $k_\mu$ on
$V^{2\lambda}$ is $\hat{\phi}^{2\lambda}_\mu$ and that of $N_{2\tau}$
on $V^{2\lambda}$ is $\hat{\theta}^{2\lambda}_{2\tau}$. 
The result follows. \eprf  

The matrix $\hat{\Theta}(2n)$ of eigenvalues of $\cB_{2n}$ 
is clearly nonsingular. Thus, 
Lemma \ref{mtrick} above shows that, given the central characters of
$S_{2n}$ and the eigenvalues of $\cB_{2n}$, and given $\mu\vdash 2n$,
we can find all the $m(\mu,2\tau), \tau\vdash n$
by solving a single system of nonsingular linear 
equations of size $p(n)\times p(n)$. 
We shall now use this
result to inductively compute the eigenvalues of 
$\cB_2,\cB_4,\ldots ,\cB_{2n}$ from the central characters 
of $S_2,S_4,\ldots ,S_{2n}$.

For $\pi\in S_{2n}$ define
$$\supp(\pi) = \{ i\in \{1,2,\ldots ,n\}\;|\; \pi(i)\not= i \mbox{ or }
                \pi(n+i)\not= n+i\;(\mbox{ or both }) \}.
$$
That is, $\supp(\pi) \cup (n + \supp(\pi))$ (here, $n + \supp(\pi) =
\{n+i\;|\;i\in \supp(\pi)\}$) is the set of end points of all
the edges of $I$ that are {\em touched} by the nontrivial 
cycles of $\pi$ (i.e., by cycles of length $\geq 2$).

Let $\mu\in \cP(2)$. For $n\geq 1$, define
\beq \label{deff}
&f(\mu,2n): \C[\cM_{2n}] \rar \C[\cM_{2n}],
\eeq
by $x\mapsto c_{\mu}(2n)\cdot x$. Note that $2n < |\mu|$ implies that
$f(\mu,2n)=0$.

Clearly $f(\mu,2n)\in \cB_{2n}$. Write
\beq \label{expf}
f(\mu,2n) &=& \sum_{\tau\in \cP(2,n)}\;
d^\tau_\mu(2n)M_{2\tau}(2n).
\eeq 

The nonnegative integers $d^\tau_\mu(2n)$ defined above can be calculated as
follows, for $n\geq |\mu|$. Below $a \vee b$ denotes the maximum 
of two nonnegative integers $a,b$. 

\bt \label{ev1}
(i) Let $\mu\in \cP(2)$ with $|\mu|=k$ and let $n\geq k$. 
For $\tau\in\cP(2,n)$ 
we have
\beqn
d^\tau_\mu(2n) &=& \left\{ \ba{cl}  
0 & \mbox{ if }|\tau| > k \\
0 & \mbox{ if }|\tau| = k \mbox{ and } \tau \not= \mu \\
2^{\ell(\mu)} & \mbox{ if }\tau=\mu 
\ea \right.
\eeqn
and, for $|\tau|=j<k$, $d^\tau_\mu(2n)$ equals
\beq \label{if}&\ds{\sum_{r=j\vee \lfloor \frac{k+1}{2} \rfloor}^{k-1}}\left\{
\ds{\sum_{s=j\vee \lfloor \frac{k+1}{2} \rfloor}^{r}}\;
(-1)^{r-s}\;\binom{r-j}{s-j}\;m((\mu,1^{2s-k}),2(\tau,1^{s-j}))\right\}
\;\binom{n-j}{r-j}&
\eeq

(ii) The set $\{ f(\mu,2n)\;|\; \mu\in \cP(2,n)\}$ is
a basis of $\cB_{2n}$.
\et
\pf 
(i) The result is clearly true if $k=0$ (in which case $f(\mu,2n)$ is the
identity map). So we may assume that $k\geq 2$. 
Let $\pi\in C_{(\mu,1^{2n-k})}$. A nontrivial $r$-cycle of $\pi$ can
touch at most $r$ edges of $I$ and thus $|\supp(\pi)|\leq k$. Moreover, if
$|\supp(\pi)|=k$ then each nontrivial $r$-cycle of $\pi$ touches exactly $r$
edges of $I$ and no edge of $I$ is touched by two distinct nontrivial
cycles. It follows that $|\supp(\pi)|=k$ implies $d(I, \pi\cdot I)=2\mu$ and
$|\supp(\pi)|\leq k-1$ implies 
$d(I,\pi\cdot I)=2(\lambda,1^{n-|\lambda|})$, 
where $\lambda\in \cP(2)$ satisfies
$|\lambda| \leq k-1$. 
Thus $d^\tau_\mu(2n)=0$ if $|\tau|>k$ or $|\tau|=k$ and $\tau\not=\mu$.

We now determine $d^\mu_\mu(2n)$. Consider the nontrivial
$r$-cycle $\sigma =(1\;2 \cdots r)\in S_{2n},\;2\leq r \leq n$. 
Then $\supp(\sigma)=\{1,2,\ldots ,r\}$ and 
$d(I,\sigma\cdot I)=2(r,1^{n-r})$. It can be checked that the only other
$r$-cycle $\pi$ with $\pi\cdot I = \sigma\cdot I$ is $\pi = (n+1\;n+r\;n+r-1
\cdots n+2)$. Since any element of $C_{(\mu,1^{2n-k})}$ has $\ell(\mu)$
nontrivial cycles it now follows from the paragraph above that
$d^\mu_\mu(2n)=2^{\ell(\mu)}$.

Now let $\tau\in\cP(2)$ with $|\tau|=j<k$. We now calculate $d^\tau_\mu(2n)$.

Fix $A\in \cM_{2n}$ with $d(I,A)=2(\tau, 1^{n-j})$ and with $I\cap A$, the intersection
of the set of edges of $I$ and $A$, given by 
$$I\cap A = \{[j+1,n+j+1], [j+2,n+j+2],\ldots ,[n,2n]\}.$$
We have
\beq \label{f1}
d^\tau_\mu(2n)&=& |\{ \pi\in C_{(\mu,1^{2n-k})}\;|\; \pi\cdot I =A\}|.
\eeq
Let $\pi\in C_{(\mu,1^{2n-k})}$ with $\pi\cdot I =A$. 
Then we clearly have
\beq \label{f2}
&\{1,2,\ldots ,j\}\seq \supp(\pi), \;\lfloor \frac{k+1}{2}\rfloor \leq
|\supp(\pi)|, \mbox{ and }|\supp(\pi)|\leq k-1,
\eeq 
where the last inequality follows from
the first paragraph of the proof.

Let $\cS(j,k,n)$ denote the set of all subsets $X$ of $\{1,2,\ldots ,n\}$
satisfying $\{1,2,\ldots ,j\}\seq X$ and 
$\lfloor \frac{k+1}{2}\rfloor \leq |X| \leq k-1$, i.e., $\cS(j,k,n)$
consists of all subsets of $\{1,2,\ldots ,n\}$ containing the elements
$\{1,2,\ldots ,j\}$ and with cardinality between $j \vee \lfloor
\frac{k+1}{2} \rfloor$ and $k-1$ (inclusive).
Partially order
$\cS(j,k,n)$ by set inclusion.

For $X\in \cS(j,k,n)$ define
\beqn
\alpha(X) &=& |\{ \pi \in C_{(\mu,1^{2n-k})}\;|\; \supp(\pi)\seq X, \;\pi\cdot
I = A\}|,\\
\beta(X) &=& |\{ \pi \in C_{(\mu,1^{2n-k})}\;|\; \supp(\pi)=X, \;\pi\cdot
I = A\}|.
\eeqn
Note that, from (\ref{f1}) and (\ref{f2}), we have
\beq \label{abeta}
d^\tau_\mu(2n)&=& \sum_{X\in \cS(j,k,n)}\beta(X).
\eeq
We have
\beqn 
\alpha(X) &=& \sum_{Y\seq X,\;Y\in \cS(j,k,n)}\beta(Y),\;\;\;X\in
\cS(j,k,n),
\eeqn
and by the principle of inclusion-exclusion
\beq \label{pie1}
\beta(X) &=& \sum_{Y\seq X,\;Y\in
\cS(j,k,n)}\;(-1)^{|X-Y|}\alpha(Y),\;\;\;X\in \cS(j,k,n).
\eeq
If $X\in \cS(j,k,n)$ with $|X|=s$, then a little reflection shows that
$$\alpha(X)=m((\mu,1^{2s-k}), 2(\tau,1^{s-j})).$$

If $X\in \cS(j,k,n)$ with $|X|=r$, then we have (from (\ref{pie1}) above)
\beq \label{pie2}
\beta(X) &=& \sum_{s=j\vee \lfloor \frac{k+1}{2}
\rfloor}^r\;(-1)^{r-s}\;\binom{r-j}{s-j}\;m((\mu,1^{2s-k}),2(\tau,1^{s-j})).
\eeq
Thus, from (\ref{abeta}) above, we have
\beqn
d^\tau_\mu(2n) &=& \sum_{X\in \cS(j,k,n)}\;\beta(X)\\ 
            &=& \sum_{r=j\vee \lfloor\frac{k+1}{2}\rfloor}^{k-1}\;\sum_{X\in
\cS(j,k,n), |X|=r}\;\beta(X).
\eeqn
Since the number of sets $X\in \cS(j,k,n)$ with $|X|=r$ is clearly
$\binom{n-j}{r-j}$ the result follows from (\ref{pie2}) above. 

(ii) This follows from the triangularity of the coefficients $d^\tau_\mu(2n)$
established in part (i) above. \eprf

Choose a linear ordering of $\cP_n$ in which the partitions are listed in
weakly increasing order of the sum of their nontrivial parts (i,e, parts
$\geq 2$). List the columns of the $\cY_n \times \cP_n$ matrix
$\hat{\Theta}(2n)$ in this order.

\bt \label{ev2}
The first column of $\hat{\Theta}(2n)$, indexed by $(1^n)$, is 
the all 1's vector.
Let $\mu\in\cP(2,n)$ with $|\mu|>0$. 
Then the column of $\hat{\Theta}(2n)$, indexed by
$(\mu,1^{n-|\mu|})$, is given by
\beqn
\left(\hat{\theta}^{2\lambda}_{2(\mu,1^{n-|\mu|})}\right)_{\lambda\vdash n} &=&
\frac{1}{2^{\ell(\mu)}}
\left\{\left(\hat{\phi}^{2\lambda}_{(\mu,1^{2n-|\mu|})}\right)_{\lambda\vdash n} - 
\ds{\sum_{\tau\in \cP(2,|\mu|-1)}} 
d^\tau_\mu(2n)\;\left(\hat{\theta}^{2\lambda}_{2(\tau,1^{n-|\tau|})}\right)_{\lambda\vdash
n}\right\} 
\eeqn
\et

\pf This follows by taking the eigenvalues on $V^{2\lambda}$ on both sides
of (\ref{expf}) and using Theorem \ref{ev1}. \eprf

\noi {\bf Proof of Theorem \ref{mt}.}
Assume the central characters of $S_2,S_4,\ldots ,S_{2n}$ and
the eigenvalues of $\cB_2, \cB_4,\ldots ,\cB_{2n-2}$ as given.

Let $\mu\in \cP(2,n)$ with $|\mu|=k$.
For $\lfloor \frac{k+1}{2}\rfloor \leq s \leq k-1$, we can, by Lemma
\ref{mtrick}, find all the nonnegative integers $m((\mu,1^{2s-k}),
2(\tau,1^{s-|\tau|})$, $\tau\in\cP(2,s)$ by solving a single system of
linear equations of size $p(s)\times p(s)$ (this requires the central
characters of $S_{2s}$ and the eigenvalues of $\cB_{2s}$ but since $s\leq
k-1\leq n-1$ the latter are known).

Thus the numbers $d^\tau_\mu(2n)$, for $\mu\in\cP(2,n),\;|\tau|<|\mu|$ can
computed from (\ref{if}).

We can now calculate the eigenvalues of $\cB_{2n}$ using the recurrence in
Theorem \ref{ev2}. 
\eprf

\bex{\em To illustrate, we calculate the eigenvalue tables $\hat{\Theta}(4)$ and
$\hat{\Theta}(6)$ starting from $\hat{\Theta}(2)$.
The
central characters of $S_4, S_6$ can be calculated from the character tables of
$S_4, S_6$ given in {\bf\cite{jk}}. 

We rewrite Lemma \ref{mtrick} as follows: for $\mu\vdash 2n$
\beq \label{rev}
(m(\mu,2\tau))_{\tau\vdash n} &=&
\hat{\Theta}(2n)^{-1}(\hat{\phi}^{2\lambda}_{\mu})_{\lambda\vdash n}.
\eeq

$\hat{\Theta}(2)$ is the $\cY_1\times \cP_1$ matrix $[1]$. Thus, from
(\ref{rev}) above we have
$$ m((2),2(1)) =  \hat{\phi}^{2(1)}_{(2)} = 1.$$
We list the elements of $\cY_2$ as$\{(2), (1,1)\}$ and the
elements of $\cP_2$ as $\{(1,1), (2)\}$. The first column of
$\hat{\Theta}(4)$ is $\left( \ba {r}1\\1 \ea\right)$.
From Theorem \ref{ev2}, the second column is 
\beqn
\left( \ba{l}
\hat{\theta}^{2(2)}_{2(2)} \\  \hat{\theta}^{2(1,1)}_{2(2)} 
\ea \right) &=& 
\frac{1}{2} \left\{ \left( \ba{l}
\hat{\phi}^{2(2)}_{(2,1^2)} \\ \hat{\phi}^{2(1,1)}_{(2,1^2)} 
\ea \right)
- d^{(0)}_{(2)}(4)\left( \ba{r} 1\\ 1 \ea\right)
\right\}. 
\eeqn
From Theorem \ref{ev1} we have
$$
d^{(0)}_{(2)}(4) = 2m((2), 2(1)) = 2.
$$
and hence the second column is $\left( \ba {r} 2\\-1 \ea\right)$ 
Thus we get
$$
\hat{\Theta}(4) = \left[ \ba {rr} 1 & 2 \\ 1 & -1 \ea \right],\;\;\;
\hat{\Theta}(4)^{-1} =  \left[ \ba {rr} 1/3 & 2/3 \\ 1/3 & -1/3 \ea \right].
$$
From (\ref{rev}) above we get
$$
\left( \ba{l} m((3,1), 2(1,1)) \\ m((3,1), 2(2)) \ea \right) =
\left[ \ba{rr} 1/3 & 2/3 \\ 1/3 & -1/3 \ea \right]
\left( \ba{l} \hat{\phi}^{2(2)}_{(3,1)}
\\ \hat{\phi}^{2(1,1)}_{(3,1)}\ea
\right) = \left( \ba{r} 0 \\ 4 \ea \right).
$$

We list the elements of $\cY_3$ as $\{(3), (2,1), (1^3)\}$ and the elements
of $\cP_3$ as $\{(1^3), (2,1), (3)\}$. The first column of $\hat{\Theta}(6)$
is the all 1's vector. 
From Theorem \ref{ev2}, the second column is
\beqn
\left( \ba{r}
\hat{\theta}^{2(3)}_{2(2,1)} \\  \hat{\theta}^{2(2,1)}_{2(2,1)}
\\  \hat{\theta}^{2(1^3)}_{2(2,1)}
\ea \right) &=& 
\frac{1}{2} \left\{ 
\left( \ba{r}
\hat{\phi}^{2(3)}_{(2,1^4)} \\ \hat{\phi}^{2(2,1)}_{(2,1^4)} 
\\ \hat{\phi}^{2(1^3)}_{(2,1^4)} 
\ea \right)
- d^{(0)}_{(2)}(6)\left( \ba{r} 1\\ 1\\1 \ea\right)
\right\}. 
\eeqn
From Theorem \ref{ev1} we have
$$
d^{(0)}_{(2)}(6) = 3m((2), 2(1)) = 3.
$$
and hence the second column is 
$\left( \ba{r} 6\\1\\-3 \ea \right)$.

From Theorem \ref{ev2}, the third  column of $\hat{\Theta}(6)$ is
\beqn
\left( \ba{l}
\hat{\theta}^{2(3)}_{2(3)} \\  \hat{\theta}^{2(2,1)}_{2(3)}
\\  \hat{\theta}^{2(1^3)}_{2(3)}
\ea \right) &=& 
\frac{1}{2} \left\{ 
\left( \ba{l}
\hat{\phi}^{2(3)}_{(3,1^3)} \\ \hat{\phi}^{2(2,1)}_{(3,1^3)} 
\\ \hat{\phi}^{2(1^3)}_{(3,1^3)} 
\ea \right)
- d^{(2)}_{(3)}(6)
\left( \ba{r}
6\\ 1 
\\ -3 
\ea \right)
- d^{(0)}_{(3)}(6)\left( \ba{r} 1\\ 1\\1 \ea\right)
\right\}. 
\eeqn
From Theorem \ref{ev1} we have
$$
d^{(0)}_{(3)}(6) = 3m((3,1), 2(1^2)) = 0,\;\;
d^{(2)}_{(3)}(6) = m((3,1), 2(2)) = 4.
$$
and hence
\beqn \hat{\Theta}(6) &=& \left[ \ba{rrr} 1&6&8\\1&1&-2\\1&-3&2 \ea
\right].
\eeqn

}\eex

We now refine the triangularity of the coefficients $d^\tau_\mu(2n)$ shown
in part (i) of Theorem \ref{ev1} above.
Define a partial order on $\cP$ as follows: $\mu\leq\lambda$ provided
$|\mu|<|\lambda|$ or $|\mu|=|\lambda|$ and $\mu$ can be obtained from
$\lambda$ by partitioning the parts of $\lambda$ into disjoint blocks and
then summing the parts in each block.
For instance,
$(5,3,2)\leq (4,2,2,1,1)$ but $(3,1,1)\not\leq (2,2,1)$ and $(2,2,1)\not\leq
(3,1,1)$.

\bl \label{pol2}
Let $\mu \in \cP(2)$ with $|\mu|=k$ and let $n\geq k$. Let $\tau\in \cP(2,n)$ 
be such that 
the coefficient $d^\tau_\mu(2n)$ defined in (\ref{expf}) above is nonzero.
Then 

\noi (i) $|\tau|\leq |\mu|$.

\noi (ii) $|\tau|=|\mu|$ implies $\tau = \mu$.

\noi (iii) $\ol{\tau} \leq \ol{\mu}$.

\el 

\pf Parts (i), (ii) follow from part (i) of Theorem \ref{ev1}.

\noi (iii) Let $\pi\in C_{(\mu,1^{2n-k})}$ with $d(I,\pi\cdot I)=2(\tau,
1^{n -|\tau|})$.
Let $\cD(I,\pi\cdot I)$ denote the (set) partition of $[2n]$ 
whose blocks are the
vertices of the connected components of the spanning subgraph of $K_{2n}$
with edge set $I\cup \pi\cdot I$ (note that each block has an even number of
elements). 
Define a (set) partition $p_{\pi}$ of $[2n]$ as follows: $i\not= j$ are in the
same block of $p_{\pi}$ provided $i=n+j$ or $j=n+i$ or $i$ and $j$ are in the same
nontrivial cycle of $\pi$. Note that each block of $p_{\pi}$ has an even
number of elements. Clearly, as set partitions, we have
\beq \label{e1}
&\cD(I,\pi\cdot I) \leq p_{\pi}.&
\eeq
Define $\mu_{\pi}$ to be the partition in $\cP(2)$ obtained from $p_{\pi}$ by taking
half the sizes of all blocks of $p_{\pi}$ of cardinality $\geq 4$. It is
easy to see, using (\ref{e1}), that
\beq \label{e2}
&|\tau|\leq |\mu_{\pi}| \leq |\mu|,& \\ \label{e3}
&|\tau|=|\mu| \mbox{ implies } \tau=\mu_{\pi}=\mu.&
\eeq
Write the parts of $\mu$ as $\{\mu_1,\ldots ,\mu_t\}$ so that the parts of
$\ol{\mu}$ are $\{\mu_1 - 1,\ldots ,\mu_t - 1\}$. 
Let $B$ be a block of $p_{\pi}$ of size $\geq 4$. Suppose this block
contains $m$ nontrivial cycles of $\pi$ whose sizes (we may assume without
loss of generality) to be $\mu_1,\ldots ,\mu_m$. Consider the hypergraph
with vertex set $B$ and edge set  the nontrivial cycles of $\pi$ contained
in $B$ together with the edges of $I$ contained in $B$. This hypergraph is
connected (since $B$ is a block of $p_\pi$) and so we have
\beqn
\frac{|B|}{2} &\leq & \mu_1 + (\mu_2 -1) + (\mu_3 - 1) +\cdots +(\mu_m - 1),
\eeqn
or, equivalently, $\frac{|B|}{2} - 1 \leq (\mu_1 -1)+\cdots + (\mu_m - 1).$

Writing the above inequality for every block of $p_{\pi}$ of size $\geq 4$
and summing we see that
\beq \label{e4}
|\ol{\mu_{\pi}}| &\leq & |\ol{\mu}|.
\eeq
The argument above also shows that  
\beq \label{e5}
|\ol{\mu_{\pi}}|=|\ol{\mu}|&\mbox{ implies }&
\ol{\mu_{\pi}}\leq\ol{\mu}.
\eeq
We now show that $ \ol{\tau} \leq \ol{\mu}$. This is clear from (\ref{e3}) if
$|\tau|=|\mu|$. Otherwise, by (\ref{e2}), $|\tau| < |\mu|$. We consider two
cases.

(a) $\cD(I,\pi\cdot I) \not= p_{\pi}$: By (\ref{e1}) and (\ref{e4}) we
have $|\ol{\tau}| < |\ol{\mu_{\pi}}| \leq |\ol{\mu}|$ and so $\ol{\tau} \leq
\ol{\mu}$.

(b) $\cD(I,\pi\cdot I) = p_{\pi}$: We have $\tau=\mu_{\pi}$.
The result follows from (\ref{e4}) and (\ref{e5}). \eprf 

We now define a polynomial in $\Q[t]$ using (\ref{if}). In Theorem \ref{eln}
below we shall evaluate this polynomial at values not covered by Theorem
\ref{ev1}. 

Given $\tau,\mu\in\cP(2)$ with $j=|\tau|\leq |\mu|=k$, define a polynomial
$\zeta^\tau_\mu(t)\in\Q[t]$ as follows:
\beqn
\zeta^\tau_\mu(t) &=& \left\{ \ba{cl}  
0 & \mbox{ if } j = k \mbox{ and } \tau \not= \mu \\
2^{\ell(\mu)} & \mbox{ if }\tau=\mu 
\ea \right.
\eeqn
and, for $j<k$, $\zeta^\tau_\mu(t)$ equals
\beqn 
&\ds{\sum_{r=j\vee \lfloor \frac{k+1}{2} \rfloor}^{k-1}}\left\{
\ds{\sum_{s=j\vee \lfloor \frac{k+1}{2} \rfloor}^{r}}\;
(-1)^{r-s}\;\binom{r-j}{s-j}\;m((\mu,1^{2s-k}),2(\tau,1^{s-j}))\right\}
\;\binom{\frac{t}{2}- j}{r-j}& 
\eeqn

\bl \label{pol1} Fix $\tau, \mu\in \cP(2)$ with $|\tau|\leq |\mu|$. Then 

\noi (i) $\zeta^{\mu}_{\mu}(t)=2^{\ell(\mu)}$.

\noi (ii) $\zeta^{\tau}_{\mu}(t)$ 
is a polynomial in $\Q[t]$ with degree $\leq |\mu|-|\tau|$. 

\noi (iii) $\zeta^\tau_\mu(t)=0$ 
unless $\ol{\tau}\leq \ol{\mu}$.

\noi (iv) $|\ol{\tau}|=|\ol{\mu}|$ implies that $\zeta^\tau_\mu(t)$ does
not depend on $t$, i.e., is a constant.

\el
\pf Parts (i) and (ii) follow from the definition of $\zeta^\tau_\mu(t)$.

\noi (iii) The result is true if $|\tau|=|\mu|$ and so we may assume
$|\tau|<|\mu|$. 
Part (iii) of Lemma \ref{pol2} and (\ref{if}) show that if
$\ol{\tau}\not\leq \ol{\mu}$ then 
$\zeta^\tau_\mu(2n)=0$ for all $n\geq |\mu|$. The result follows. 

\noi (iv) The result is true if $|\tau|=|\mu|$ and so we may assume
$|\tau|<|\mu|$. Let $|\ol{\tau}|=|\ol{\mu}|$ and let $n\geq |\mu|$.
Let $\pi,\sigma \in C_{(\mu,1^{2n-k})}$
satisfy $d(I,\pi\cdot I)= 2(\tau, 1^{n-|\tau|})$ 
and $\sigma\cdot I = \pi\cdot I$. Then, by (\ref{e1}) and by case (a) in
the proof of part (iii) in Lemma \ref{pol2} above,
we have $\cD(I,\sigma\cdot I)= \cD(I,\pi\cdot I) = p_{\pi}$. 
It follows that $\zeta^\tau_\mu(2n)$ does not depend on $n$. The result follows. \eprf

\bt \label{eln} Let $\mu\in\cP(2)$ with $|\mu|=k$.

\noi (i) For $k\leq n$ we have
\beqn
f(\mu,2n)&=&2^{\ell(\mu)}M_{2\mu}(2n) + 
\sum_{\tau\in\cP(2,k-1)}\;\zeta^\tau_\mu(2n)M_{2\tau}(2n).
\eeqn

\noi (ii) For $n < k \leq 2n$ we have
\beqn
f(\mu,2n)&=&\sum_{\tau\in\cP(2,k-1)}\;\zeta^\tau_\mu(2n)M_{2\tau}(2n).
\eeqn

\noi (iii) For $2n < k$ and $\tau\in\cP(2),\;|\tau|\leq n$,
we have $\zeta^\tau_\mu(2n)=0$.
\et
\pf (i) This follows from Theorem \ref{ev1}.

Before proving parts (ii) and (iii) we make the following observation.

Let $2n\geq k$ so that $c_\mu(2n)\not= 0$. Then, just as in the proof of part (i)
of Theorem \ref{ev1} we have
\beq \label{ex}
f(\mu,2n)&=& d^\mu_\mu(2n)M_{2\mu}(2n) +
\sum_{\tau\in\cP(2,k-1)}d^\tau_\mu(2n)M_{2\tau}(2n).
\eeq
Fix $\tau\in\cP(2,k-1)$ with $|\tau|=j<k$. Define $\gamma_\tau$ to be the
number of perfect matchings $A$ in $\cM_{2j}$ with 
$d(\{[1, j+1], [2, j+2], \ldots ,[j, 2j]\}, A)=2\tau$. Thus the number
of perfect matchings in $\cM_{2n}$ with $d(I,A) = 2(\tau, 1^{n-j})$ is
$\gamma_\tau \binom{n}{j}$.

For $j\vee \lfloor \frac{k+1}{2} \rfloor \leq r \leq k$ define
\beqn
\alpha(r,\tau)&=&\{\pi\in C_{(\mu, 1^{2n-k})}\;|\;\mbox{supp}(\pi)=\{1,2,\ldots
,r\},\;d(I,\pi\cdot I)=2(\tau,1^{r-j})\}.
\eeqn
A little reflection shows that
\beq \nonumber
d^{\tau}_{\mu}(2n) &=& 
\frac{
{\ds{\sum_{r=j\vee \lfloor \frac{k+1}{2} \rfloor}^k}} \alpha(r,\tau){\binom{n}{r}}
}
{
\gamma_\tau{\binom{n}{j}}
}\\ \label{ex1}
&=&
{\ds{\sum_{r=j\vee \lfloor \frac{k+1}{2} \rfloor}^k}} 
\;\frac{j!}{r!}\;\frac{\alpha(r,\tau)}{\gamma_\tau}\;(n-j)(n-j-1)\cdots
(n-r+1).
\eeq
The expression in (\ref{ex1}) above is valid for all $n\geq k/2$ and thus it
follows that
\beq \label{ex2}
\zeta^\tau_\mu(t)&=&
{\ds{\sum_{r=j\vee \lfloor \frac{k+1}{2} \rfloor}^k}} 
\;\frac{j!}{r!}\;\frac{\alpha(r,\tau)}{\gamma_\tau}\;\left(\frac{t}{2}-j\right)
\left(\frac{t}{2}-j-1\right)\cdots \left(\frac{t}{2}-r+1\right).
\eeq

(ii) Since $n<k$ we clearly have $d^\mu_\mu(2n)=0$. The result now follows
from (\ref{ex}), (\ref{ex1}), and (\ref{ex2})  above.

(iii) This follows from (\ref{ex2}) on noting that, for $2n<k$ and
$|\tau|=j\leq n$ we have $n\in \{j,j+1,\ldots , 
\lfloor \frac{k+1}{2} \rfloor - 1\}$.
\eprf

{{\bf  \section {Content evaluation of symmetric functions}}}

We now consider algorithms for expressing
$\phi^\lambda_{\mu}$ and 
$\theta^{2\lambda}_{2\mu}$, for fixed $\mu\in \cP(2)$ and 
varying $\lambda\in \cY$, as content evaluations of symmetric functions.
The motivation comes from certain basic
results in the representation theory of symmetric groups {\bf\cite{cst2,g,ov}}.
We now recall these in items (I)-(III) below (this will also be used
in the next section on eigenvectors).

\begin{enumerate}
\item[(I)] Consider an irreducible $S_n$-module $V^{\lambda}$, for
$\lambda\in \cY_n$.
Since the branching is multiplicity free, the decomposition into irreducible
$S_{n-1}$-modules of $V^{\lambda}$ is
canonical.
Each of these modules, in turn, decompose canonically into
irreducible $S_{n-2}$-modules. Iterating this construction, we get a
canonical decomposition of $V^{\lambda}$ into irreducible $S_1$-modules,
i.e., one
dimensional subspaces. Thus, there is a canonical basis of $V^{\lambda}$,
determined
up to scalars, and  called the
{\em Gelfand-Tsetlin (or GZ-) basis} of $V^{\lambda}$.

\item[(II)] For $i=1,2,\ldots ,n$ define $X_i = (1, i) + (2, i) + \cdots +
(i-1, i) \in \C[S_n]$. The $X_i$'s are called the
{\em Young-Jucys-Murphy elements (YJM-elements)}. Note that $X_1=0$.

Consider
the Fourier transform, i.e., the algebra isomorphism
\beq\label{iso}
\C[S_n]&\cong&\bigoplus_{\lambda\in \cY_n} \mbox{End}(V^{\lambda}),
\eeq
given by
$$\pi \mapsto ( V^{\lambda} \buildrel {\pi}\over \rightarrow
V^{\lambda}\;:\;
\lambda
\in \cY_n),\;\;\pi\in S_n.$$
We have identified a canonical basis, the GZ-basis, in each
$S_n$-irreducible. Let
$\mbox{D}(V^{\lambda})$ consist of all operators on $V^{\lambda}$
diagonal in the GZ-basis of $V^{\lambda}$. It is known that the image of
$\bigoplus_{\lambda\in \cY_n} \mbox{D}(V^{\lambda})$ (a maximal commutative
subalgebra of the right hand side of (\ref{iso})) under the inverse  
Fourier transform is the subalgebra of $\C[S_n]$ generated by $X_1,\ldots
,X_n$, which is thus a maximal commutative subalgebra of $\C[S_n]$. 
It follows that the only common eigenvectors of $X_1,\ldots ,X_n$ in
an irreducible $V^{\lambda}$ are (up to scalars) 
the elements of the GZ-basis of
$V^{\lambda}$. Moreover, the eigenvalues of the YJM elements    
on the GZ-basis vectors
in each irreducible can also be written down once we parametrize the
GZ-basis by standard Young tableaux. We recall this in the next item below.

\item[(III)] Let $\mu\in \cY$. A {\em Young tableau of shape $\mu$}
is obtained by taking the Young diagram $\mu$ and filling its
$|\mu|$ boxes 
(bijectively) with the numbers $1,2,\ldots ,| \mu |$.
A Young tableau  is said to be {\em standard} if the numbers in
the boxes strictly increase along each row and each column of the Young
diagram of $\mu$. Let $\tab(n,\mu)$, where $\mu \in
\cY_n$,
denote the set of all
standard Young tableaux of shape $\mu$ and let $\tab(n)=\cup_{\mu\in
\cY_n}\tab(n,\mu)$. There is a well known bijection between
$\tab(n,\lambda)$ and sequences $(\lambda_1,\lambda_2,\ldots ,\lambda_n)$
of Young diagrams with $\lambda_n = \lambda$ and $\lambda_{i}\in
\lambda_{i+1}^-$, for $1\leq
i \leq n-1$ (given $T\in \tab(n,\lambda)$, define $\lambda_i$ to be the 
diagram obtained by considering the boxes of $T$ containing the numbers 
$1,\ldots ,i$). It now easily follows from the branching rule that 
the GZ-basis of $V^{\lambda}$ can be parametrized by $\tab(n,\lambda)$.
Given $T\in \tab(n,\lambda)$, we write $v_T$ for the corresponding GZ-basis
vector of $V^{\lambda}$.

Given $T\in \tab(n,\lambda)$, the
eigenvalue of $X_i$ on $v_T$ is $c(b_T(i))$, the content of 
the box $b_T(i)$ of $T$ containing $i$. 

\end{enumerate}

Let $f=f(X_1,\ldots ,X_n)$ be a symmetric polynomial in $X_1,\ldots ,X_n$.
By considering the $GZ$-basis of $V^{\lambda}$ we see that the action of $f$
on $V^{\lambda}$ is multiplication by the scalar $f(c(\lambda))$. Using the
Fourier transform, it now follows that any symmetric polynomial in
$X_1,\ldots ,X_n$ is in $Z[\C[S_n]]$. The converse
of this assertion is also true.

Given $n$ variables $x_1,\ldots ,x_n$ and $1\leq k \leq n$, we let
$e_k(x_1,\ldots ,x_n)$ denote the elementary symmetric polynomials. Suppose
$a=\{a_1,\ldots , a_n\}$ and $b=\{b_1,\ldots ,b_n\}$ are two multisets of
(complex) numbers of cardinality $n$. By considering the polynomials
$(x-a_1)\cdots (x-a_n)$ and $(x-b_1)\cdots (x-b_n)$ we see that $a=b$ as
multisets if and only if $e_k(a_1,\ldots ,a_n)=e_k(b_1,\ldots ,b_n)$, for
$1\leq k \leq n$.

Let $\lambda,\mu \in\cY_n$. The number of 0's
in $c(\lambda)$ is the number of boxes in the main diagonal of $\lambda$,
the number of 1's is the number of boxes in the first superdiagonal, the
number of -1's is the number of boxes in the first subdiagonal and so on. It
follows that $\mu = \lambda$ if and only if $c(\mu)=c(\lambda)$ if and only
if $e_k(c(\lambda))=e_k(c(\mu))$ for $1\leq k \leq n$.

Fix $\lambda \in \cY_n$. For $1\leq k \leq n$, define the following
symmetric polynomials in $X_1,\ldots ,X_n$:
\beqn f_k(X_1,\ldots ,X_n) &=& \prod_{\mu} (e_k(X_1,\ldots ,X_n) -
e_k(c(\mu)),
\eeqn
where the product is over all $\mu\in \cY_n$ with
$e_k(c(\mu))\not=e_k(c(\lambda))$.

Let $\mu\in \cY_n$ with $\mu\not= \lambda$. Then, by the observation above,
$e_k(c(\mu))\not=e_k(c(\lambda))$ for some $1\leq k \leq n$. It follows that
\beqn
\left(\prod_{k=1}^nf_k(X_1,\ldots , X_n)\right)\cdot V^{\mu} &=&
\left\{ \ba{ll}  0 & \mbox{ if $\mu\in \cY_n,\;\mu\not=\lambda$,} \\
                 \mbox{nonzero scalar} & \mbox{ if $\mu=\lambda$.}\ea
\right.
\eeqn
Using the Fourier transform we see that 
every element in $Z[\C[S_n]]$ is a
symmetric polynomial in $X_1,\ldots ,X_n$. 

Thus $Z[\C[S_n]]$ consists of all symmetric polynomials 
in $X_1,\ldots ,X_n$. This is Jucys' fundamental 
theorem given, with a different proof, in {\bf\cite{j}}. Constructive proofs 
of this result are given in 
Murphy {\bf\cite{mur}}, Moran {\bf\cite{mo}}, 
Diaconis and Greene {\bf\cite{dg}}, and Garsia {\bf\cite{g}}. A good
reference for this material is the book of
Cecchereni-Silberstein, Scarabotti, and Tolli {\bf\cite{cst2}}.

In the
present context we have the following basic problem: for fixed
$\mu\in \cP(2)$, write the conjugacy class
sum $c_\mu(n)\in Z[\C[S_n]]$ as a linear combination of, say, the
power sum symmetric functions in $X_1,\ldots ,X_n$ and say something about
the dependence of the coefficients on $n$.

Given $f\in \Lambda[t]$ and $n\geq 1$ we define the {\em YJM
evaluation} $f(n,X)$ to be the element of $Z[\C[S_n]]$ obtained from $f$ by
setting $t=n,\;x_i=0\mbox{ for }i>n$, and $x_i=X_i,\;i=1,\ldots ,n$.

The following result was proved in {\bf\cite{cgs}}. An algorithm
for constructing the symmetric function $W_\mu$ was given in {\bf\cite{g}}.
See {\bf\cite{cst2}} for another proof ( Part (iv) below is taken from
Theorem 5.4.7 of this reference).

\bt \label{garsia}
For each $\mu\in\cP(2)$ there is an algorithm to compute a symmetric
function $W_{\mu}\in \Lambda[t]$ such that

\noi (i) $\{W_{\mu}\;:\;\mu\in\cP(2)\}$ is a $\Q[t]$-module basis of $\Lambda$.

\noi (ii) For $\mu\in\cP(2)$ and  $n\geq 1$ we have
\beqn
W_{\mu}(n,X) &=& c_\mu(n).
\eeqn

\noi (iii) For $\mu\in\cP(2)$ and $\lambda\in \cP$ we have
\beqn
W_{\mu}(c(\lambda)) &=& \phi^{\lambda}_{\mu}.
\eeqn

\noi (iv) Let $\mu\in\cP(2)$ with multiplicity of $i$ equal to $m_i$, $i\geq 2$.
The expansion of $W_{\mu}$ in the power sum basis has the form
\beqn
W_{\mu} &=& \sum_{\lambda\;\leq\;\ol{\mu}}\;a^{\lambda}_{\mu}(t)\;p_{\lambda},
\eeqn
where

(a) $a^{\lambda}_{\mu}(t)\in \Q[t]$ with degree $\leq
\frac{|\ol{\mu}|-|\lambda|}{2} + \ell(\ol{\mu}) - \ell(\lambda)$.

(b) $a_{\mu}^{\ol{\mu}}=\frac{1}{\prod_{i\geq 2} m_i!}$ and
$a_{\mu}^{\lambda}\in\Q$ (i.e., does not depend on $t$) for
$|\lambda|=|\ol{\mu}|$.

(c) $a^\lambda_\mu(t)=0$ if $|\ol{\mu}|$ and $|\lambda|$ do not have the
same parity.

\et

\noi {\bf Remark} Let $\mu,\tau\in\cP(2)$. Using Theorem \ref{garsia}(i),
we can write
\beqn
W_\mu W_\tau &=& \sum_\lambda \omega^\lambda_{\mu,\tau}(t) \;W_\lambda,
\eeqn
where the sum is over finitely many $\lambda\in \cP(2)$ and
$\omega^\lambda_{\mu,\tau}(t) \in \Q[t]$. From Theorem
\ref{garsia}(iii) we have
\beqn
c_\mu(n) c_\tau(n) &=& \sum_\lambda \omega^\lambda_{\mu,\tau}(n) \;
c_\lambda(n),\;\;n\geq 1.
\eeqn
In other words, the structure constants of the algebra of fixed conjugacy
classes (the so-called {\em Farahat-Higman algebra})
are integer valued rational polynomials. See {\bf\cite{cgs,cst2}} for more
details.

\bt \label{mks1}
For each $\mu\in\cP(2)$ there is an algorithm to compute a symmetric
function $E_{\mu}\in \Lambda[t]$ such that

\noi (i) $\{E_{\mu}\;:\;\mu\in\cP(2)\}$ is a $\Q[t]$-module basis of
$\Lambda[t]$.

\noi (ii) For $\mu\in\cP(2)$ and $\lambda\in \cP$ we have
\beqn
E_{\mu}(c(2\lambda)) &=& \theta^{2\lambda}_{2\mu}.
\eeqn

\noi (iii) Let $\mu\in\cP(2)$ with multiplicity of $i$ equal to $m_i$,
$i\geq 2$.
The expansion of $E_{\mu}$ in the power sum basis has the form
\beqn
E_{\mu} &=& \sum_{\lambda\;\leq\;\ol{\mu}}\;b^{\lambda}_{\mu}(t)\;p_{\lambda},
\eeqn
where

(a) $b^{\lambda}_{\mu}(t)\in \Q[t]$ with degree $\leq
|\ol{\mu}|-|\lambda| + \ell(\ol{\mu}) - \ell(\lambda)$.

(b) $b_{\mu}^{\ol{\mu}}=\frac{1}{2^{\ell(\mu)}\prod_{i\geq 2} m_i!}$ and
$b_{\mu}^{\lambda}\in\Q$ (i.e., does not depend on $t$) for
$|\lambda|=|\ol{\mu}|$.

\et

\noi {\bf Remark} Let $\mu,\tau\in\cP(2)$. Using Theorem \ref{mks1}(i),
we can write
\beqn
E_\mu E_\tau &=& \sum_\lambda \beta^\lambda_{\mu,\tau}(t) \;E_\lambda,
\eeqn
where the sum is over finitely many $\lambda\in \cP(2)$ and
$\beta^\lambda_{\mu,\tau}(t) \in \Q[t]$. From Theorem
\ref{mks1}(ii) we have
\beqn
M_{2\mu}(2n) M_{2\tau}(2n) &=& \sum_\lambda \beta^\lambda_{\mu,\tau}(2n) \;
M_{2\lambda}(2n),\;\;n\geq 1.
\eeqn
In other words, the structure constants of the algebra of fixed orbitals
are integer valued rational polynomials. For a direct study of these
structure constants
in more detail see the two recent papers {\bf\cite{ac,co,t}} (our focus 
in this paper is more on
the eigenvalues and eigenvectors of $\cB_{2n}$).
These papers work in
the context of the Hecke algebra of the Gelfand pair $(S_{2n},H_n)$ (which
explains the extra factor $2^n n!$ in their structure constants). 

\pf We proceed by induction on $|\mu|$. 
Set $E_{(0)}=1$ and
assume that, for some $k\geq 1$, we have defined $E_{\mu} \in \Lambda[t]$, 
for all $\mu\in \cP(2)$
with $|\mu|\leq k-1$, such that 
items (ii) and (iii)
in the statement of the theorem are satified.

Now let $\mu\in \cP(2)$ with $|\mu|=k$ and with multiplicity of $i$ equal to
$m_i,\;i\geq 2$. 
Define
\beqn
E_{\mu} &=& \frac{1}{2^{\ell(\mu)}}\;\left\{W_{\mu} -
\sum_{\tau\in\cP(2,k-1)}
\zeta^{\tau}_{\mu}(t)\;E_{\tau}\right\}.
\eeqn
We shall now verify items (ii) and (iii)(a), (iii)(b) 
in the statement for $E_{\mu}$. 
We begin with item (iii).

By Theorem \ref{garsia}(iv) we can write 
\beq \label{t1}
W_{\mu} &=& \sum_{\lambda\;\leq\;\ol{\mu}}\;a^{\lambda}_{\mu}(t)\;p_{\lambda},
\eeq
where degree of $a^{\lambda}_{\mu}(t) 
\leq
\frac{|\ol{\mu}|-|\lambda|}{2} + \ell(\ol{\mu}) - \ell(\lambda)
\leq
|\ol{\mu}|-|\lambda| + \ell(\ol{\mu}) - \ell(\lambda)$.

Let $\tau\in \cP(2)$ with $|\tau|\leq k-1$. 
By the induction hypothesis we can write
\beq  \label{t2}
E_{\tau} &=& \sum_{\lambda\;\leq\;\ol{\tau}}\;b^{\lambda}_{\tau}(t)
\;p_{\lambda},
\eeq
where degree of $b^{\lambda}_{\tau}(t) \leq
|\ol{\tau}|-|\lambda| + \ell(\ol{\tau}) - \ell(\lambda)$.

Now, degree of $\zeta^\tau_\mu(t) \leq |\mu| - |\tau| =
|\mu|-|\ol{\tau}| - \ell(\ol{\tau})$ (by Lemma \ref{pol1}(ii)) 
and thus degree of
$\zeta^\tau_\mu(t)b^\lambda_\tau(t)$$ \leq |\mu| - |\lambda| - \ell(\lambda)
= |\ol{\mu}|-|\lambda|+\ell(\ol{\mu})-\ell(\lambda)$.

By Lemma \ref{pol1}(iii) we have that $\zeta^\tau_\mu(t)\not= 0$ implies 
$\ol{\tau}\leq \ol{\mu}$. Item  (iii)(a)
now follows from (\ref{t1}) and (\ref{t2}). 

Item  (iii)(b) also
follows from (\ref{t1}) and (\ref{t2}) by using the induction hypothesis,
Theorem \ref{garsia} (iv)(b)  and
Lemma \ref{pol1}(iii).

We now verify item (ii). Let $\lambda\in \cY_m$ and consider the following
three cases:

\noi (a) $k\leq m$: This follows from Theorems \ref{garsia}(iii), Theorem
\ref{eln}(i),
and the induction hypothesis.

\noi (b) $m < k \leq 2m$: 
We need to show that $E_\mu(c(2\lambda))=0$.
This follows from
Theorem \ref{garsia}(iii) and Theorem \ref{eln}(ii).

\noi (c) $k>2m$: We need to show that $E_\mu(c(2\lambda))=0$.
By Theorem \ref{garsia}(iii) we have $W_\mu(c(2\lambda))=0$. By
the induction hypothesis $E_\tau(c(2\lambda))=0$ for $m<|\tau|$ and by 
Lemma \ref{eln}(iii) $\zeta^\tau_\mu(2m)=0$ for $|\tau|\leq m$. The result
follows.

That completes the proof of items (ii) and (iii). Item (i) now follows from
Theorem \ref{garsia}(i) and the triangular definition of the $E_\mu$. 
\eprf

It is easily seen that property (ii) of Theorem \ref{mks1} characterizes the
symmetric function $E_\mu$.

\bco Let $f,g\in \Lambda[t]$. Suppose that there exists $n_0$ 
such that $f(c(2\lambda))=g(c(2\lambda))$ for all
$\lambda\in \cY,\;|\lambda|\geq n_0$. Then $f=g$. 
\eco

\pf Suppose $f\not= g$. Write
\beq \label{k1}
f-g &=& a_{\mu_1}(t)E_{\mu_1}+ a_{\mu_2}(t)E_{\mu_2}+\cdots +
a_{\mu_k}(t)E_{\mu_k},
\eeq
where $\mu_i\in \cP(2)$ for all $i$ and $a_{\mu_i}(t)\not= 0$ for all $i$.

Choose a positive integer $m$ such that $m\geq n_0$, $|\mu_i|\leq m$ for
$1\leq i \leq k$, and $a_{\mu_i}(2m)\not= 0$ for $1\leq i \leq k$.
We can now rewrite (\ref{k1}) (by adding terms with zero coefficients) as
\beq \label{k2}
f-g& = & \sum_{\mu\in\cP(2,m)}\;a_{\mu}(t)E_{\mu},
\eeq
where not all $a_{\mu}(2m)$ are zero.

Evaluate both sides of (\ref{k2}) on the contents of $2\lambda$, for every
$\lambda\vdash m$. By assumption we get
\beq \label{k3}
0& = & \sum_{\mu\in\cP(2,m)}
\;a_{\mu}(2m)\hat{\theta}^{2\lambda}_{2(\mu,1^{m-|\mu|})},\;\;\lambda\vdash m.
\eeq
From (\ref{k3}) we get that a nontrivial linear combination of
the columns of (the nonsingular matrix) $\hat{\Theta}(2m)$ is zero, 
a contradiction. \eprf

\bex\label{wep} {\em 
Below we give tables of $W_{\mu}$ and $E_{\mu}$ polynomials for $|\mu|\leq
4$. The $W_{\mu}$ polynomials are from {\bf\cite{cgs,g}} while the $E_{\mu}$
polynomials were calculated by hand using Theorem \ref{mks1}.

\beqn 
W_0 = 1 &&
E_0 = 1 \\\\
W_2 = p_1&&
E_2 = \frac{p_1}{2}- \frac{t}{4}\\\\
W_3 = p_2 -\frac{t(t-1)}{2} &&
E_3 = \frac{p_2}{2} - p_1 + \frac{3t-t^2}{4}\\\\
W_{2,2} = \frac{p_1^2}{2} - \frac{3p_2}{2} + \frac{t(t-1)}{2}&&
E_{2,2} = \frac{p_1^2}{8} - \frac{3p_2}{4} + \frac{(10-t)p_1}{8} +
\frac{9t^2-24t}{32}\\\\
W_4= p_3 - (2t-3)p_1 && E_4 = \frac{p_3}{2} -\frac{9p_2}{4} 
+ \frac{(11-2t)p_1}{2} + \frac{8t^2 - 23t}{8}
\eeqn

We can calculate the eigenvalue table $\hat{\Theta}(8)$ of $\cB_8$ using the
list above. We list the elements of $\cP_4$ in the order $\{(1^4), (2,1^2),
(2,2), (3,1), (4)\}$ and the elements of $\cY_4$ in the order $\{(4), (3,1),
(2,2), (2,1^2), (1^4)\}$. We have
\beqn \hat{\Theta}(8) &=& \left[ \ba{rrrrr} 
1&12&12&32&48\\
1&5&-2&4&-8\\
1&2&7&-8&-2\\
1&-1&-2&-2&4\\
1&-6&3&8&-6 
\ea
\right].
\eeqn

}\eex

{{\bf  \section {  Eigenvectors: similar algorithms for
$\hat{\phi}^\lambda_\mu$ and $\hat{\theta}^{2\lambda}_{2\mu}$}}}

In this section we shall give an inductive procedure to write down a specific 
eigenvector (a so called first GZ-vector) in each eigenspace of the (left)
actions of $Z[\C[S_n]]$ and $\cB_{2n}$ on $\C[S_n]$ and $\C[\cM_{2n}]$
respectively. This then yields simple inductive algorithms to calculate
$\hat{\phi}^\lambda_\mu$ and $\hat{\theta}^{2\lambda}_{2\mu}$ (that do not
depend on knowing the symmetric group characters).

To begin with, it will be useful to know (as suggested
by (\ref{rrs}) and Lemma \ref{mfree}) how
GZ-vectors behave under restriction and induction.

The case of restriction follows from the following result.

\bl \label{res} Let $\lambda\in \cY_n$ and consider the irreducible
$S_n$-module $V^\lambda$.

(i) Let $v\in V^\lambda$ be an eigenvector for the action of $X_1,\ldots
,X_{n-1}$. The $v$ is also an eigenvector for the action of $X_n$.

(ii)  Suppose $T\in \tab(n,\lambda)$ and $v\in V^\lambda$ satisfy
$$X_i\cdot v = c(b_T(i))v,\;1\leq i \leq n-1.$$
Then $X_n\cdot v = c(b_T(n))v$.

(iii) The GZ-basis of $V^{\lambda}$ is the union
of the GZ-bases of $V^{\mu}$, as $\mu$ varies over $\lambda^-$.

\el
\pf (i) Let $X$ be the sum of all transpositions in $S_n$. Note that
$X=X_1+\cdots +X_n$ and that $X$ is in the center of $\C[S_n]$. Thus, by
Schur's lemma, 
the action of $X$ on $V^\lambda$ is multiplication by a scalar. Thus
$v$ is an eigenvector for the action of $X_n = X - (X_1+\cdots + X_{n-1})$.

(ii) The action of $X$ on $V^{\lambda}$ is multiplication by a scalar
$\alpha$. By considering a GZ-vector of $V^{\lambda}$ we see that $\alpha$
is equal to the 
sum of the contents of all boxes of the Young diagram $\lambda$.
The result follows.

(iii) This follows from parts (i) and (ii) above using the branching rule
(\ref{brres}).
\eprf

Now we consider the case of induction. Since we will also be applying this
construction to the case of the regular module $\C[S_n]$, which is not
multiplicity free, we first extend the notion of a GZ-vector to a
$S_n$-module with a single isotypical component.

Let $V$ be a $S_n$-module with a single isotypical component, the
irreducibles occuring in $V$ all being isomorphic to $V^\lambda,\;\lambda\in
\cY_n$. Let $T\in \tab(n,\lambda)$ and define the following subspace of $V$:
$$V_T = \{v\in V\;|\;X_i(v)=c(b_T(i))v,\;i=1,\ldots ,n\}.$$
It is easy to see that we have the canonical decomposition:
$$ V = \oplus_{T\in \tab(n,\lambda)} V_T.$$
By a { \em GZ-vector of $V$ associated to $T$} we mean a
nonzero vector in $V_T$.

For a Young diagram $\lambda$ let 
$\cO(\lambda)$ be the set of boxes corresponding to the outer
corners of $\lambda$. Note that no two boxes in $\cO(\lambda)$ have the same
content. For $\lambda\in \cY_n$, we 
denote the isotypical component of $V^\lambda$ in a $S_n$-module $W$ by
$W^\lambda$.

\bl \label{ind} Let $W$ be a $S_n$-module and let
$$U = \C[S_{n+1}] \otimes_{\C[S_n]} W
=\ind_{\;S_n}^{\;S_{n+1}}(W).$$

Let $T\in\tab(n,\lambda)$ and let $v\in W^\lambda$ be a GZ-vector associated
to $T$.
Let $\mu\in \lambda^+$ and let $b\in \cO(\lambda)$ be
the box added to $\lambda$ to get $\mu$. Let $S\in \tab(n+1,\mu)$ be the
standard tableau obtained from $T$ by adding $n+1$ in box $b$. 

Then 
$$ \prod_{d\in \cO(\lambda)\setminus\{b\}}(X_{n+1} - c(d)I)\cdot
(1\otimes v)$$
is a GZ-vector of $U^\mu$ associated to $S$.
\el
\pf It suffices to prove the case $W=V^\lambda$. In this case $v=v_T$ (upto
scalars) and, by the branching rule, $U=\oplus_{\tau\in \lambda^+} V^\tau$. 
Clearly, $1\otimes v_T \in U$ is $\not= 0$. Write 
$1\otimes v_T=\sum_{\tau\in \lambda^+} v_\tau$, where $v_\tau \in V^\tau$.

For $1\leq i \leq n$ we have $X_i\cdot (1\otimes v_T) = 
1\otimes (X_i\cdot v_T) = c(b_T(i))(1\otimes v_T)$. It follows that
$X_i \cdot v_\tau = c(b_T(i))v_\tau$, $1\leq i \leq n,\;\tau\in \lambda^+$. 
From part (ii)
of Lemma \ref{res} it now follows that $X_{n+1}\cdot v_\tau = c(d)
v_\tau,\;\tau\in \lambda^+$,
where $d$ is the box added to $\lambda$ to get $\tau$. 
The result follows. \eprf

Let $\lambda = (\lambda_1,\ldots ,\lambda_t)\vdash n$. Define the standard
tableau $R\in \tab(n,\lambda)$ by filling the boxes of $\lambda$ with the
integers $1,2,\ldots ,n$ in {\em row major order}, i.e., the first row is
filled with the numbers $1,2,\ldots , \lambda_1$ (from left to right), the
second row with the numbers $\lambda_1+1,\lambda_1+2,\ldots ,\lambda_2$ and
so on. Given a $S_n$-module $W$, a nonzero vector $v$ in 
$(W^\lambda)_R$ will be called  
{\em a first Gelfand-Tsetlin vector} in $W^\lambda$.

We now give an example of a first GZ-vector and rederive a result from
{\bf\cite{gm1,l}}. First, we make a definition. The {\em perfect matching
derangement operator}
$$D_{2n}:\C[\cM_{2n}] \rar \C[\cM_{2n}]$$
is defined as follows: for $A\in \cM_{2n}$ set $D_{2n}(A) = \sum_B B$, where
the sum is over all $B\in \cM_{2n}$ with $d(A,B)$ having no part equal to 2.
In other words, $D_{2n}=\sum_\mu N_{2\mu}$, where the sum is over all
$\mu\in \cP(2)$ with $|\mu|=n$. For $\lambda\vdash n$, let
$m^{2\lambda}_{2n}$ denote the eigenvalue of $D_{2n}$ on $V^{2\lambda}$.

Fix a matching $A\in \cM_{2n}$. The number of $B\in \cM_{2n}$ with $d(A,B)$
having no part equal to 2 is easily seen (by inclusion-exclusion) to be
$$d(2n)=\sum_{i=0}^{n-1}(-1)^i \binom{n}{i} (2n-2i-1)!!$$

We denote by $v_{2\lambda}$ the first GZ-vector in the subspace
$V^{2\lambda}$ of $\C[\cM_{2n}]$. For the rest of this section
fix $J=\{[1,2],[3,4],\ldots ,[2n-1,2n]\}\in \cM_{2n}$.

\bex \label{glm} {\em (i) Clearly, $v_{2(n)} = \sum_{A\in \cM_{2n}} A$.
Let $\mu\in\cP_n$.
The coefficient of $J$ in $v_{2(n)}$ is $1$ while the coefficient of $J$ 
in $N_{2\mu}(v_{2(n)})$ (respectively, $D_{2n}(v_{2(n)})$) is 
$|\cM(J,2\mu)|$ (respectively, $d(2n)$). It follows that
\beqn
\hat{\theta}^{2(n)}_{2\mu} &=& |\cM(J,2\mu)|,\\
m^{2(n)}_{2n} &=& d(2n).
\eeqn
It is easy to see that
$$|\cM(J,2\mu)| = \frac{2^n n!}{z_{\mu}2^{\ell(\mu)}}.$$

(ii) We now write down $v_{2(n-1,1)}$. Using the inductive structure of
$\C[\cM_{2n}]$ given in Lemma \ref{mfree} (v), (vi) and applying 
Lemmas \ref{res} and \ref{ind}, we get from item (i) above,
\beqn
v_{2(n-1,1)} &=& \left(X_{2n-1} - (2n-2)I\right)\cdot\left(
\sum_{A\in \cM_{2n-2}} (A\cup
\{[2n-1,2n]\})\right)\\
  &=& \left( \sum_{i=1}^{2n-2} \sum_{A\in \cM_{2n}, [i,2n]\in A} A\right)
  - (2n-2)\left(\sum_{A\in \cM_{2n}, [2n-1,2n]\in A} A\right)\\
  &=& \sum_{A\in \cM_{2n}} A - (2n-1)
                \left(\sum_{A\in \cM_{2n}, [2n-1,2n]\in A} A\right).
\eeqn
The coefficient of $J$ in $v_{2(n-1,1)}$ is $-(2n-2)$ and the coefficient of
$J$ in $D_{2n}(v_{2(n-1,1)})$ is $d(2n)$. We can easily
calculate the coefficient of $J$ in $N_{2\mu}(v_{2(n-1,1)})$. Two cases
arise:

(a) $\mu$ has no part equal to $1$: 
The coefficient of $J$ in $N_{2\mu}(v_{2(n-1,1)})$ 
is $|\cM(J,2\mu)|$.

(b) $\mu$ has a part equal to $1$: Let $\mu'\in \cP_{n-1}$ be obtained from
$\mu$ by deleting a $1$ from the parts of $\mu$ and let
$J'=\{[1,2],[3,4],\ldots ,[2n-3,2n-2]\}\in \cM_{2n-2}$. 
The coefficient of $J$ in $N_{2\mu}(v_{2(n-1,1)})$ is 
$|\cM(J,2\mu)|- (2n-1)|\cM(J',2\mu')|$.

It follows that 
\beqn
\hat{\theta}^{2(n-1,1)}_{2\mu} &=& \left\{ \ba{ll}
\frac{|\cM(J,2\mu)|}{-(2n-2)}, & \mbox{ if $1$ is not a part of $\mu$,}\\
&\\
\frac{|\cM(J,2\mu)| - (2n-1)|\cM(J',2\mu')|}{-(2n-2)}, 
& \mbox{ if $1$ is a part of $\mu$}.
\ea \right.
\eeqn
and that $$m^{2(n-1,1)}_{2n} = \frac{d(2n)}{-(2n-2)}.$$
}\eex

In principle, it is possible to extend the method of Example \ref{glm} to
certain other eigenspaces, such as $V^{2(n-2,2)}$ and $V^{2(n-2,1,1)}$, and
derive complicated explicit formulas for $m^{2(n-2,2)}_{2n}$ and
$m^{2(n-2,1,1)}_{2n}$. We do not pursue this here. 
Instead we shall show how Lemmas \ref{res} and
\ref{ind} can be used to give a practical
recursive
algorithm for calculating $\hat{\theta}^{2\lambda}_{2\mu}$. To be efficient
we shall not write down the eigenvectors explicitly 
but only keep track of the values of
these eigenvectors at a (subexponential) number of linear functionals on
$\C[\cM_{2n}]$.

For $\mu\in \cP_n$ define a linear functional
$$f_{2\mu} : \C[\cM_{2n}] \rar \C$$
as follows: given $v\in \C[\cM_{2n}]$ write
$$v=\sum_{A\in \cM_{2n}} \alpha_A A,\;\;\alpha_A\in \C.$$
Define $f_{2\mu}(v) = \sum_A \alpha_A,$ where  the sum is over all $A\in
\cM_{2n}$ with $d(J,A)=2\mu$. 
We call $(f_{2\mu}(v))_{\mu\vdash n}$ the {\em
orbital coefficients} of $v\in \C[\cM_{2n}]$. Note that the vector $v$,
living in a vector space of dimension $(2n-1)!!$, has
only $p(n)$ orbital coefficients.

Given $\lambda\in \cY_n$, 
let $v_{2\lambda}$ denote the first GZ-vector of
the submodule $V^{2\lambda}$ of $\C[\cM_{2n}]$, normalized so that the
coefficient of $J$ in $v_{2\lambda}$ is $1$. Then it follows that
$$\hat{\theta}^{2\lambda}_{2\mu} = f_{2\mu}(v_{2\lambda}).$$

Thus, the eigenvalues can be determined once we know the orbital
coefficients of the first GZ-vectors. The basic idea of the algorithm is
to inductively compute the orbital
coefficients using Lemmas \ref{res} and \ref{ind}. This leads  to the following problem,
called the {\em update problem}:

Given the orbital coefficients of $v\in \C[\cM_{2n}]$, determine the orbital
coefficients of $X_n\cdot v$.

In order to solve the update problem we need to go slightly beyond orbital
coefficients to relative orbital coefficients.  

Let 
$$\cP_n' = \{(\mu,i)\;|\;\mu\in \cP_n\mbox{ and $i$ is a part of $\mu$}\}.$$
Elements of $\cP_n'$ are called {\em pointed partitions} of $n$. 
Let $pp(n)$ denote
the number of pointed partitions of $n$. Clearly, $pp(n) = 1 +p(1)+\cdots
+p(n-1)$ (note that $pp(n)$ is also subexponential. We have $pp(13)=272$). 
Pointed partitions play an important role in Okounkov-Vershik
theory (see {\bf\cite{ov,cst2}})  as $pp(n)$ is the dimension of the relative
commutant $\{\pi\in \C[S_n]\;|\;\pi \C[S_{n-1}] = \C[S_{n-1}] \pi \}$.

For $(\mu,i)\in \cP_n'$ define a linear functional
$$f_{(2\mu, 2i)} : \C[\cM_{2n}] \rar \C$$
as follows: given $v\in \C[\cM_{2n}]$ write
$$v=\sum_{A\in \cM_{2n}} \alpha_A A,\;\;\alpha_A\in \C.$$
Define $f_{(2\mu,2i)}(v) = \sum_A \alpha_A,$ where the sum is over all $A\in
\cM_{2n}$ with $d(J,A)=2\mu$ and with the size of the component of $J\cup A$
containing the edge $[2n-1,2n]$ being $2i$. 
We call $(f_{(2\mu, 2i)}(v))_{(\mu,i)\in \cP_n'}$ the {\em
relative orbital coefficients} of $v\in \C[\cM_{2n}]$. 

For $\lambda\in \cY_n, \mu\in \cP_n$ we now have 
$$\hat{\theta}^{2\lambda}_{2\mu} = \sum_i f_{(2\mu,2i)}(v_{2\lambda}),$$
where the sum is over all parts $i$ of $\mu$.

The update problem for relative orbital coefficients can be easily  solved
using the following lemma.
\bl \label{upl} 
Let $A\in\cM_{2n}$. 
Let $C_1,C_2,\ldots ,C_t$ be the components of the spanning
subgraph of $K_{2n}$ with edge set $J\cup A$, with $C_t$ containing the edge
$[2n-1,2n]$. Let $2\mu_i$ be the number of vertices of $C_i$, $i=1,\ldots
,t$. Thus $\{2\mu_1,\ldots ,2\mu_t\}$ is the multiset of parts of $d(J,A)$.

(i) Let $s$ be a vertex of $C_j$, $j=1,\ldots , t-1$ and put $A'=
(s,2n-1)\cdot A$. Then the multiset of
parts of $d(A', J)$ is 
$$(\{2\mu_1,\ldots ,2\mu_t\} - \{2\mu_j , 2\mu_t\}) \cup \{2(\mu_j + \mu_t)\},$$ 
with $2(\mu_j+\mu_t)$ as
the size of the component of $A' \cup J$ containing the edge
$[2n-1,2n]$.

(ii) Traverse the vertices of the alternating cycle $C_t$ in cyclic order,
beginning at the vertex $2n$ and going towards $2n-1$. List the
vertices encountered as $\{2n, 2n-1, i_1, i_2, \ldots ,i_{2k-1}, i_{2k}\}$,
where $k\geq 0$ and $2\mu_t = 2k+2$. Then

(a) Let $j\in \{1,2,\ldots ,k\}$ and put $A' = (i_{2j}, 2n-1)\cdot A$.
The multiset of parts of $d(A', J)$ is 
$\{2\mu_1,\ldots ,2\mu_{t-1},2\mu_t -2j, 2j\}$, with
$2\mu_t - 2j$ as the size of the component of $A'\cup J$ containing the edge
$[2n-1, 2n]$.

(b) Let $j\in \{1,2,\ldots ,k\}$ and put $A' = (i_{2j-1}, 2n-1)\cdot A$.
The multiset of parts of $d(A', J)$ is $\{2\mu_1,\ldots ,2\mu_t\}$, with
$\mu_t$ as the size of the component of $A'\cup J$ containing the edge
$[2n-1, 2n]$.

\el
\pf (i) Let $[s,x], [2n-1,y]\in A$. Then $A' = (A\setminus \{[s,x],
[2n-1,y]\}) \cup \{[2n-1,x], [y,s]\}$. It follows that 
$C_k$, $k\in \{1,\ldots ,t-1\}\setminus \{j\}$, continue to remain
components of $J\cup A'$ and that $C_j$ and $C_t$ merge into a single
alternating cycle in $J\cup A'$.

(ii)(a) It is clear that $C_1,\ldots ,C_{t-1}$ continue to be components of
$J\cup A'$ and that $C_t$ splits into two alternating cycles with vertex
sets 
$$\{i_{2j+1}, i_{2j+2}, \ldots ,i_{2k-1}, i_{2k}, 2n-1, 2n\} \mbox{ and } 
\{i_1, i_2, \ldots ,i_{2j-1}, i_{2j}\}.$$

(ii)(b) Similar to case (ii)(a) except that $C_t$ does not split. $\Box$

For $v\in \C[\cM_{2n}]$, define 
\beqn [v] &=& (f_{(2\mu,2i)}(v))_{(\mu,i)\in\cP_n'}
\eeqn
to be the vector of the relative orbital coefficients of $v$. We denote
$f_{(2\mu,2i)}(v)$ by $v(2\mu,2i)$.

The following is the algorithm for updating the vector
of relative orbital coefficients. Its correctness directly follows from
Lemma \ref{upl}. 

{\texttt
{\bf{Algorithm 1.}} {\em (Update for relative orbital coefficients)}

{\bf{Input}}: $[v]$, for some $v\in \C[\cM_{2n}]$, and an integer $a$.

{\bf{Output}}: $[u]$, where $u = (X_{2n-1} - aI )\cdot(v) \in \C[\cM_{2n}]$.
        We denote the output by $F_{a}([v])$.

{\bf{Method}}:

1. For all $(\mu,i)\in \cP_n'$ do $\gamma(2\mu,2i) = 0$.

2. For all $(\mu,i)\in \cP_n'$ do

 $\;\;\;\;\;\;\;$2a. Write the multiset of parts of $\mu$ as
$\{\mu_1,\mu_2,\ldots ,\mu_t\}$, where $\mu_t = i$.  

 $\;\;\;\;\;\;\;$2b. For $j = 1 \mbox{ to } t-1$ do

 $\;\;\;\;\;\;\;\;\;\;\;\;\;\;$2b.1. 
$\mu' = 
(\{\mu_1, \mu_2,\ldots ,\mu_t\} \setminus \{\mu_j ,\mu_t\})
\cup \{\mu_j + \mu_t\},\;\;i'=\mu_j + \mu_t$.

$\;\;\;\;\;\;\;\;\;\;\;\;\;\;$2b.2. 
$\gamma(2\mu',2i') = 2\mu_j v(2\mu,2i) + \gamma(2\mu', 2i').$

 $\;\;\;\;\;\;\;$2c. $k = \mu_t - 1$.

 $\;\;\;\;\;\;\;$2d. For $j = 1 \mbox{ to } k$ do

$\;\;\;\;\;\;\;\;\;\;\;\;\;\;$2d.1. 
$\mu' = 
(\{\mu_1, \mu_2,\ldots ,\mu_{t-1},\mu_t -j, j\}, 
\;\;i'=\mu_t - j$.

$\;\;\;\;\;\;\;\;\;\;\;\;\;\;$2d.2. 
$\gamma(2\mu',2i') = v(2\mu,2i) + \gamma(2\mu', 2i').$

 $\;\;\;\;\;\;\;\;\;\;\;\;\;\;$2d.3. 
$\gamma(2\mu,2i) = v(2\mu,2i) + \gamma(2\mu, 2i).$

3. For all $(\mu,i)\in \cP_n'$ do $u(2\mu,2i) = \gamma(2\mu,2i) - a v(2\mu,2i)$.

4. RETURN $(u(2\mu,2i))_{(\mu,i)\in \cP_n'}$.

}

We now give the inductive algorithm for computing the rows of the eigenvalue
tables of $\hat{\Theta}(2n)$. In Step 5 below we use the convention that,
for a proposition $P$, $[P]$ equal to 1 if $P$ is true and
is equal to 0 if $P$ is false.

{\bf{Algorithm 2}}

{\bf Input:} (i) $\lambda'\in \cY_{n+1}$, 
with $\lambda = \lambda' -\mbox{ \{last box in
last row of $\lambda$\}}\in \cY_n$.

(ii) The row of $\hat{\Theta}(2n)$ indexed by $\lambda$, i.e., 
$(\hat{\theta}^{2\lambda}_{2\mu})_{\mu\in \cP_n}$.

{\bf Output:} The row of $\hat{\Theta}(2n+2)$ indexed by $\lambda'$, i.e.,
$(\hat{\theta}^{2\lambda'}_{2\mu'})_{\mu'\in \cP_{n+1}}$. 

{\bf Method:} 

1. For all $(\mu',i)\in \cP_{n+1}'$ do $v(2\mu',2i) = 0$. 
   
2. For all $\mu\in \cP_n$ do $v(2\mu\cup \{2\},2) = \hat{\theta}^{2\lambda}_{2\mu}$.

3. Let the Young diagram $2\lambda$ have $k+1$ outer boxes. 
   Adding two boxes (in a row) in the\\ 
   \hspace*{0.43in}place of one of these outer boxes 
   yields $2\lambda'$. Denote the $k$ other outer boxes by
   $b_1,\ldots ,b_k$. 

4. For $j = 1 \mbox{ to } k$ do $[v] = F_{c(b_j)}([v])$  

5. For all $\mu'\in \cP_{n+1}$ do $\hat{\theta}^{2\lambda'}_{2\mu'} =
\frac{\ds{\sum_{i=1}^{n+1}}\;[i\mbox{ is a part of $\mu'$} ]
\;v(2\mu',2i)}{v(2(1^{n+1}),2)}$  

6. RETURN $(\hat{\theta}^{2\lambda'}_{2\mu'})_{\mu'\in\cP_{n+1}}$

\bt Algorithm 2 is correct.
\et
\pf Let $u\in \C[\cM_{2n}]$ be the first GZ-vector in $V^{2\lambda}$,
normalized so that the coefficient of $J$ is 1. Let $v\in \C[\cM_{2n+2}]$ be
the vector corresponding to $1\otimes u \in
\ind^{S_{2n+1}}_{S_{2n}}(\C[\cM_{2n}])$, under the isomorphism between
$\ind^{S_{2n+1}}_{S_{2n}}(\C[\cM_{2n}])$ and
$\res^{S_{2n+2}}_{S_{2n+1}}(\C[\cM_{2n+2}])$ (Lemma \ref{mfree}(v,vi)).
Then it follows that steps 1 and 2 of Algorithm 2 correctly calculate $[v]$.

It now follows from Lemma \ref{ind} that steps 3, 4, 5, and 6 of Algorithm 2
correctly compute the orbital coefficients of the first GZ-vector of
$V^{2\lambda'}$. \eprf

We have implemented Algorithms 1 and 2 in Maple. Both the program and its
binary file are available at {\bf\cite{sr}}. The program is able to compute
$\hat{\theta}^{2\lambda}_{2\mu}$ reasonably quickly for $|\lambda|=|\mu|\leq
20$. We were able to determine the entire spectrum of $D_{40}$.

\bex {\em
We give below the eigenvalue table $\hat{\Theta}(10)$ computed using this
program.
We list the elements of $\cP_5$ in the order $\{(1^5), (2,1^3),
(2^2,1), (3,1^2), (3,2), (4,1), (5)\}$ and the elements of $\cY_5$ 
in the order $\{(5), (4,1), (3,2), (3,1^2),
(2^2,1), (2,1^3), (1^5)\}$. We have
\beqn \hat{\Theta}(10) &=& \left[ \ba{rrrrrrr} 
1&20&60&80&160&240&384\\
1&11&6&26&-20&24&-48\\
1&6&11&-4&20&-26&-8\\
1&3&-10&2&-4&-8&16\\
1&0&5&-10&-10&10&4\\
1&-4&-3&2&10&6&-12\\
1&-10&15&20&-20&-30&24
\ea
\right].
\eeqn
Summing the fifth and seventh columns of $\hat{\Theta}(10)$ we get the
spectrum  of $D_{10}$:
\beqn
&m^{(10)}_{10} = 544,
m^{(8,2)}_{10} = -68,
m^{(6,4)}_{10} = 12,
m^{(6,2,2)}_{10} = 12,&\\
&m^{(4,4,2)}_{10} = -6,
m^{(4,2,2,2)}_{10} = -2,
m^{(2,2,2,2,2)}_{10} = 4.&
\eeqn
Note the sole zero value in row 5, column 2. The eigenvalue table
$\hat{\Theta}(2n)$ tends to have far fewer zero values than the character
(or central character) table of $S_n$. For instance, $p(15)=176$ and of the
$176^2=30976$ entries in the character table of $S_{15}$ as many as $11216$
are zero while only $878$ of the entries in $\hat{\Theta}(30)$ are zero.
}\eex

Recently, Ku and Wong {\bf\cite{kw}} gave elegant explicit formulas for
$m^{2(1^n)}_{2n}$ and $m^{2(2^m,1^{n-2m})}_{2n}$. Namely, they showed that
$$m^{2(1^n)}_{2n}=(-1)^{n-1}(n-1),\;\;m^{2(2^m,
1^{n-2m})}_{2n}=(-1)^{n-2}((m-1)n - m^2 + 2m +1).$$ 
It would be interesting to see whether these formulas can be derived from the
algorithm presented here. This possibility arises as follows. The number $k$ 
of
times the for loop in Step 4 of Algorithm 2 is executed depends on the
number of outer boxes of the input Young diagram. In the case of the Young
diagrams $2(1^n)$ and $2(2^m,1^{n-2m})$ this number is 1 or 2 throughout
(i.e., at every level of recursion). This considerably simplifies the
recursion and it may be possible to use generating function techniques to
derive the formulas above. We hope to return to this later.

We shall now give an almost identical algorithm for computing the central
characters of $S_n$, based on the inductive structure
(\ref{rrs}) of the regular modules $\C[S_n]$.

For $\mu\in \cP_n$ define a linear functional
$$g_{\mu} : \C[S_n] \rar \C$$
as follows: given $v\in \C[S_n]$ write
$$v=\sum_{\pi\in S_n} \alpha_\pi \pi,\;\;\alpha_\pi\in \C.$$
Define $g_{\mu}(v) = \sum_\pi \alpha_\pi,$ where  the sum is over 
$\pi\in C_\mu$. 
We call $(g_{\mu}(v))_{\mu\vdash n}$ the {\em
class coefficients} of $v\in \C[S_n]$. Note that the vector $v$,
living in a vector space of dimension $n!$, has
only $p(n)$ class coefficients.

Given $\lambda\in \cY_n$, 
let $v_{\lambda}$ denote a first GZ-vector in
the submodule $(\C[S_n])^\lambda$ of $\C[S_n]$, normalized so that the
coefficient of the identity permutation in $v_{\lambda}$ is $1$. 
Then it follows that
$$\hat{\phi}^{\lambda}_{\mu} = g_{\mu}(v_{\lambda}).$$

Thus, the eigenvalues can be determined once we know the class
coefficients of first GZ-vectors. The basic idea of the algorithm is
to inductively compute the class
coefficients using Lemma \ref{ind}. Like before, this leads  
to the {\em update problem}:

Given the class coefficients of $v\in \C[S_n]$, determine the class
coefficients of $X_n\cdot v$.

To solve the update problem we define relative class coefficients.  
For $(\mu,i)\in \cP_n'$ define a linear functional
$$g_{(\mu, i)} : \C[S_n] \rar \C$$
as follows: given $v\in \C[S_n]$ write
$$v=\sum_{\pi\in S_n} \alpha_\pi \pi,\;\;\alpha_\pi\in \C.$$
Define $g_{(\mu,i)}(v) = \sum_\pi \alpha_\pi$, where
the sum is over all $\pi\in S_n$ 
with $\pi\in C_\mu$ and with the size of the cycle of $\pi$
containing $n$ being $i$. 
We call $(g_{(\mu, i)}(v))_{(\mu,i)\in \cP_n'}$ the {\em
relative class coefficients} of $v\in \C[S_n]$. 

For $\lambda\in \cY_n, \mu \in \cP_n$ we now have 
$$\hat{\phi}^{\lambda}_{\mu} = \sum_i g_{(\mu,i)}(v_{\lambda}),$$
where the sum is over all parts $i$ of $\mu$.

The update problem for relative class coefficients can be easily solved
using the following lemma.
\bl \label{uplc} 
Let $\pi\in S_n$ with 
$C_1,C_2,\ldots ,C_t$ as its disjoint cycles and with
$C_t$ containing $n$. Let $\mu_i = |C_i|,\;i=1,\ldots ,t$, so that
$\{\mu_1,\ldots ,\mu_t\}$ is the multiset of cycle lengths of $\pi$.

(i) Let $s$ be an element of $C_j$, $j=1,\ldots , t-1$ and put $\pi'=
(s,n) \pi$. Then the multiset of cycle lengths of $\pi'$ is
$$(\{\mu_1,\ldots ,\mu_t\} - \{\mu_j , \mu_t\}) \cup \{\mu_j + \mu_t\},$$ 
with $\mu_j+\mu_t$ as the length of the cycle containing $n$.

(ii) Write $C_t = (n, i_k, i_{k-1},\ldots ,i_1)$, where $k\geq 0$ and $\mu_t=k+1$. 
   Let $j\in \{1,2,\ldots ,k\}$ and put $\pi' = (i_j, n)\cdot \pi$.
Then the  multiset of parts of $\pi'$ is 
$\{\mu_1,\ldots ,\mu_{t-1},\mu_t -j, j\}$, with
$\mu_t - j$ as the length of the cycle containing $n$.

\el
\pf This is similar to the proof of Lemma \ref{upl}. $\Box$

For $v\in \C[S_n]$, define 
\beqn [v] &=& (g_{(\mu,i)}(v))_{(\mu,i)\in\cP_n'}
\eeqn
to be the vector of the relative class coefficients of $v$. We denote
$g_{(\mu,i)}(v)$ by $v(\mu,i)$.

The following is the algorithm for updating the vector
of relative class coefficients. Its correctness directly follows from
Lemma \ref{uplc}. 

{\texttt
{\bf{Algorithm 3.}} {\em (Update for relative class coefficients)}

{\bf{Input}}: $[v]$, for some $v\in \C[S_n]$, and an integer $a$.

{\bf{Output}}: $[u]$, where $u = (X_{n} - aI )\cdot(v) \in \C[S_n]$.
        We denote the output by $G_{a}([v])$.

{\bf{Method}}:

1. For all $(\mu,i)\in \cP_n'$ do $\gamma(\mu,i) = 0$.

2. For all $(\mu,i)\in \cP_n'$ do

 $\;\;\;\;\;\;\;$2a. Write the multiset of parts of $\mu$ as
$\{\mu_1,\mu_2,\ldots ,\mu_t\}$, where $\mu_t = i$.  

 $\;\;\;\;\;\;\;$2b. For $j = 1 \mbox{ to } t-1$ do

 $\;\;\;\;\;\;\;\;\;\;\;\;\;\;$2b.1. 
$\mu' = 
(\{\mu_1, \mu_2,\ldots ,\mu_t\} \setminus \{\mu_j ,\mu_t\})
\cup \{\mu_j + \mu_t\},\;\;i'=\mu_j + \mu_t$.

$\;\;\;\;\;\;\;\;\;\;\;\;\;\;$2b.2. 
$\gamma(\mu',i') = \mu_j v(\mu,i) + \gamma(\mu', i').$

 $\;\;\;\;\;\;\;$2c. $k = \mu_t - 1$.

 $\;\;\;\;\;\;\;$2d. For $j = 1 \mbox{ to } k$ do

$\;\;\;\;\;\;\;\;\;\;\;\;\;\;$2d.1. 
$\mu' = 
(\{\mu_1, \mu_2,\ldots ,\mu_{t-1},\mu_t -j, j\}, 
\;\;i'=\mu_t - j$.

$\;\;\;\;\;\;\;\;\;\;\;\;\;\;$2d.2. 
$\gamma(\mu',i') = v(\mu,i) + \gamma(\mu', i').$

3. For all $(\mu,i)\in \cP_n'$ do $u(\mu,i) = \gamma(\mu,i) - a v(\mu,i)$.

4. RETURN $(u(\mu,i))_{(\mu,i)\in \cP_n'}$.
}

We now give the inductive algorithm for computing the rows of the central
character tables of $S_n$.

{\bf{Algorithm 4}}

{\bf Input:} (i) $\lambda'\in \cY_{n+1}$, 
with $\lambda = \lambda' -\mbox{ \{last box in
last row of $\lambda$\}}\in \cY_n$.

(ii) The row of the central character table of $S_n$ indexed by $\lambda$, i.e.,
$(\hat{\phi}^{\lambda}_{\mu})_{\mu\in \cP_n}$.

{\bf Output:} The row of the central character table of $S_n$ 
indexed by $\lambda'$, i.e.,
$(\hat{\phi}^{\lambda'}_{\mu'})_{\mu'\in \cP_{n+1}}$. 

{\bf Method:} 

1. For all $(\mu',i)\in \cP_{n+1}'$ do $v(\mu',i) = 0$. 

2. For all $\mu\in \cP_n$ do $v(\mu\cup \{1\},1) = 
\hat{\phi}^{\lambda}_{\mu}$.

3. Let $\lambda$ have $k+1$ outer boxes. One of these outer boxes, when
   added to $\lambda$, yields $\lambda'$. \\
\hspace*{0.43in}Denote the $k$ other outer boxes by
   $b_1,\ldots ,b_k$. 

4. For $j = 1 \mbox{ to } k$ do $[v] = G_{c(b_j)}([v])$  

5. For all $\mu'\in \cP_{n+1}$ do
$\hat{\phi}^{\lambda'}_{\mu'} =
\frac{\ds{\sum_{i=1}^{n+1}}\;[i\mbox{ is a part of $\mu'$} ]\;v(\mu',i)}
{v(1^{n+1},1)}$  

6. RETURN $(\hat{\phi}^{\lambda'}_{\mu'})_{\mu'\in\cP_{n+1}}$

\bt Algorithm 4 is correct.
\et 
\pf This is similar to the proof of correctness of Algorithm 2.
\eprf

This algorithm has also been implemented in {\bf\cite{sr}}.

\end{document}